\documentclass[11pt]{article}
\usepackage{mathrsfs}
\usepackage{amssymb}
\usepackage{amsfonts}
\usepackage{color,xcolor}
\usepackage{graphicx}
\usepackage{manfnt}
\usepackage{graphicx}
\newtheorem{theorem}{Theorem}[section]
\newtheorem{corollary}{Corollary}[section]
\newtheorem{lemma}{Lemma}[section]

\newtheorem{definition}{Definition}[section]
\newtheorem{remark}{Remark}[section]
\newtheorem{example}{Example}[section]

\def\[{{\Big[}}\def\]{{\Big]}}\def\<{{\langle}}\def\>{{\rangle}}\def\({{\Big(}}
\def\){{\Big)}}

\def\sgn{\mbox{\rm sgn}}

\def\={&\!\!=\!\!&}
\def\cC{{\mathcal C}}\def\cD{{\mathcal D}}

\def\cL{{\mathcal L}}
\def\cM{{\mathcal M}}

\def\mN{{\mathbb N}}
\def\mR{{\mathbb R}}

\def\geq{\geqslant}\def\leq{\leqslant}
\def\div{\mathord{{\rm div}}}

\begin{document}
\title{\bf Kinetic solutions \\ for nonlocal scalar conservation laws}
\author{Jinlong Wei$^a$, Jinqiao Duan$^b$and Guangying Lv$^c$\thanks{Corresponding author Email: gylvmaths@henu.edu.cn} }
\date{$^a$ School of Statistics and Mathematics, Zhongnan University of
\\ Economics and Law, Wuhan, Hubei 430073, China \\ $^b$ Department of Applied Mathematics\\ Illinois Institute of Technology, Chicago, IL 60616
\\ $^c$School of Mathematics and Statistics, Henan University
\\
Kaifeng, Henan 475001, China}
 \maketitle
\noindent{\hrulefill}
\vskip1mm\noindent
{\bf Abstract} This work is devoted to examine the uniqueness and existence of kinetic solutions for a class of scalar conservation laws involving a nonlocal super-critical diffusion operator. Our proof for uniqueness is based upon the analysis on a   microscopic contraction functional  and the existence is enabled by a parabolic approximation. As an illustration, we obtain the existence and uniqueness of kinetic solutions for the generalized fractional Burgers-Fisher type equations. Moreover, we demonstrate the kinetic solutions' Lipschitz continuity   in time, and   continuous dependence on   nonlinearities and L\'{e}vy measures.

   \vskip1mm\noindent
{\bf Keywords:}  Kinetic solution; Nonlocal conservation laws; Uniqueness; Existence; Anomalous diffusion

  \vskip2mm\noindent
{\bf MSC (2010):} 35L03; 35L65; 35R11
 \vskip0mm\noindent{\hrulefill}

\section{Introduction}\label{sec1}\setcounter{equation}{0}

The present paper is concerned with the  anomalous  diffusion related to  the L\'{e}vy flights \cite{Str,SZF, Duan2015}.  At the macroscopic modeling level, this means the   Laplacian for normal diffusion is replaced by a fractional power of the (negative) Laplacian.  We consider the following  partial differential equation, coupling   a conservation law with an anomalous  diffusion:
\begin{eqnarray}\label{1.1}
\frac{\partial}{\partial t} \rho(t,x)+\div_xF(\rho)+\nu\cL \rho=0, \ \ (t,x)\in (0,T) \times \mR^d,
\end{eqnarray}
fulfilling the initial data
\begin{eqnarray}\label{1.2}
\rho(t=0,x)=\rho_0(x), \ \ x\in \mR^d,
\end{eqnarray}
where $\nu$ is a nonnegative parameter and
\begin{eqnarray}\label{1.3}
\rho_0\in L^1\cap BV(\mR^d), \ F\in W^{1,\infty}_{loc}(\mR;\mR^d), \ \cL=(-\Delta_x)^{\frac{\alpha}{2}} \ \mbox{or} \ \sum_{i=1}^d(-\partial^2_{x_i,x_i})^{\frac{\alpha}{2}}, \ \alpha \in (0, 1).
\end{eqnarray}
Moreover, $(-\Delta_x)^{\frac{\alpha}{2}}$ is the nonlocal or fractional
Laplacian in $\mR^d$ (see \cite{DI}), defined, for any $\varphi \in \cD(\mR^d)$, $x\in \mR^d$, by
\begin{eqnarray}\label{1.4}
(-\Delta_x)^{\frac{\alpha}{2}}\varphi(x)=c(d,\alpha) \mbox{P.V.} \int_{\mR^d}\frac{\varphi(x)-\varphi(z+x)}
{|z|^{d+\alpha}}dz=c(d,\alpha)\int_{\mR^d}\frac{\varphi(x)-\varphi(z+x)}
{|z|^{d+\alpha}}dz,
\end{eqnarray}
with  $c(d,\alpha)=\alpha2^{\alpha-1}\pi^{-d/2}\Gamma{(\frac{d+\alpha}{2})}/\Gamma{(\frac{2-\alpha}{2})}$.
We also denote  $(-\partial^2_{x_i,x_i})^{\frac{\alpha}{2}}$ for the 1-dimensional fractional
Laplacian.

\vskip1mm \par

The nonlocal Cauchy problem (\ref{1.1})-(\ref{1.2}) has attracted a lot of attention  for the past few years due to its broad  applications in
mathematical finance \cite{DI}, hydrodynamics \cite{SK}, acoustics \cite{Sug}, trapping effects in surface diffusion \cite{ZA}, statistical mechanics \cite{Za,BW}, relaxation phenomena \cite{SW}, physiology \cite{BFW,DGV} and molecular biology \cite{Dro,Ali}, and its relation with stochastic analysis  \cite{BKW1,BKW2,Zha}.

We briefly mention some    recent works on well-posedness of (\ref{1.1})-(\ref{1.2}), which are   relevant for the present paper.
We first recall a remarkable result on the scalar conservation law without diffusion ($\nu=0$):
\begin{eqnarray}\label{1.5}
\frac{\partial}{\partial
t} \rho(t,x)+\mbox{div}_xF(\rho)=0, \ \ (t,x)\in (0,T) \times \mR^d.
\end{eqnarray}
Since (\ref{1.5}) is hyperbolic, classical solutions, starting out from smooth initial values, spontaneously develop discontinuities. Hence, in general, only weak solutions may exist. But weak solutions may fail to be unique in general.   By introducing an entropy formulation
\begin{eqnarray}\label{1.6}
\frac{\partial}{\partial t}\eta(\rho)+\div_xQ(\rho)\leq 0,
\end{eqnarray}
Kru\u{z}kov  \cite{Kru} showed the uniqueness results for entropy solutions in  $L^\infty$ space.
 \vskip1mm\par
The general Kru\u{z}kov type theory on well-posedness for nonlocal version of (\ref{1.5}), i.e. (\ref{1.1}) with $\alpha\in (1,2)$ (called sub-critical) was initiated by \cite{BFW} for the fractional Burgers equation ($\nu>0$ and $F(\rho)=a\rho^r$) in
Bessel potential and/or Morrey spaces. This result was then strengthened   by Droniou, Gallou\"{e}t and Vovelle \cite{DGV};  using a splitting method, they proved the global existence and uniqueness of regular solutions. A general result in this
direction was obtained by Droniou and Imbert   \cite{DI}, by means of the  ``reverse maximum principle" and Duhamel's formula; they proved  the
existence and uniqueness for regular solutions to the
Hamilton-Jacobi equation.
 \vskip1mm\par
The critical ($\alpha=1$) and super-critical ($\alpha\in (0,1)$) cases are more difficult.  Alibaud \cite{Ali}  obtained  well-posedness results for $L^\infty$-solutions of fractional conservation laws.
\vskip1mm\par
Recently,  Lions, Perthame and Tadmor \cite{LPT} proved that, if $\rho$ is an entropy solution and belongs to $L^1$ space, then for any $v\in \mR$, $u(t,x,v)$ defined by
\begin{eqnarray}\label{1.7}
u(t,x,v)=\chi_{\rho}(v)=1_{(0,\rho(t,x))}(v) -1_{(\rho(t,x),0)}(v)
\end{eqnarray}
satisfies
\begin{eqnarray}\label{1.8}
\frac{\partial}{\partial t}u(t,x,v)+f(v)\cdot \nabla_xu(t,x,v)=\frac{\partial}{\partial v}m(t,x,v), \ \ (t,x,v)\in (0,T) \times \mR^d\times \mR,
\end{eqnarray}
in $\cD^\prime((0,T)\times\mR^{d+1})$ and initial data
\begin{eqnarray}\label{1.9}
u(t=0)=\chi_{\rho_0}(v), \ \ (x,v)\in \mR^d\times \mR,
\end{eqnarray}
where $f=F^\prime$ and $m$ is a nonnegative measure.
But when discussing (\ref{1.8})-(\ref{1.9}), $L^1$ is a natural space  for the solutions. Based upon this observation, Perthame extended Kru\u{z}kov $L^\infty$ theory for entropy solutions and developed an $L^1$ theory for kinetic solutions (\cite{Per,Per1}).  How to generalize this $L^1$ theory to the Cauchy problem (\ref{1.1})-(\ref{1.2})  is an interesting issue.
    \vskip1mm\par
As claimed in \cite{Ali}, one can define "intermediate" (for  $"classical\Rightarrow entropy \Rightarrow intermediate \Rightarrow weak"$) solutions for (\ref{1.1}) by
\begin{eqnarray}\label{1.10}
\frac{\partial}{\partial t}\eta(\rho)+\div_xQ(\rho)
+\nu(-\Delta_x)^{\frac{\alpha}{2}}\eta(\rho)\leq 0.
\end{eqnarray}
  Previously mentioned works did not use this entropy formulation, since the doubling variable technique is not appropriate to this solution,  and to a very great degree, intermediate solution is non-unique. Furthermore, as inspired by \cite{Per,Per1},  we note that (\ref{1.10}) may be  suitable for us to establish a relationship between (\ref{1.1}) and the following nonlocal linear convection-diffusion equation
\begin{eqnarray}\label{1.11}
\frac{\partial}{\partial t}u(t,x,v)+f(v)\cdot \nabla_xu+ \nu\cL u=\frac{\partial}{\partial v}(m+n), \ \ (t,x,v)\in (0,T) \times \mR^d\times \mR,
\end{eqnarray}
via a kinetic formulation, with certain nonnegative measures $m$ and $n$.  When we deal with (\ref{1.11}), some technical difficulties may be overcome in order to show the uniqueness for kinetic solutions.
On account of this fact, in the present paper we introduce a notion of kinetic solution (analogue of \cite{Per}) and will prove that under the assumption (\ref{1.3}), the Cauchy problem (\ref{1.1})-(\ref{1.2}) is well-posed. It is non-trivial to get
the uniqueness of the kinetic solution to (\ref{1.1})-(\ref{1.2}) because of the nonlocal term $\cL \rho$,
see Section 3. Moreover, we revisit the continuous dependence on nonlinearities and L\'{e}vy measures. Comparing with
the results in \cite{ACJ1,ACJ2}, we delete the assumption $\rho_0\in L^\infty$.
   \vskip1mm\par
This paper is organized as follows. In Section
2, we introduce some notions on solutions for (\ref{1.1})-(\ref{1.2}), and then prove the uniqueness and existence of  kinetic solutions in Section 3. We further discuss
the regularity properties and continuous dependence (on nonlinearities and L\'{e}vy measures) for kinetic solutions in Section 4.

\section{Entropy solutions and kinetic solutions}\label{sec2}
\setcounter{equation}{0}
We take $\nu>0$ and the analysis on $\nu\cL$ is the same as $\nu(-\Delta_x)^{\frac{\alpha}{2}}$, for writing simplicity, we choose $\cL=(-\Delta_x)^{\frac{\alpha}{2}}$ in the present paper,
 and we take $\nu=1$ in Section 2 and Section 3.  Now we introduce some notions.

\begin{definition}  (Entropy solution) Let (\ref{1.3}) hold and $\rho_0(x)\in L^\infty(\mR^d)$.  A function $\rho\in L^\infty([0,T)\times \mR^d)\cap \cC([0,T];L^1(\mR^d))\cap L^\infty([0,T);BV(\mR^d))$ is said to be an entropy solution of (\ref{1.1})-(\ref{1.2}), if for every smooth convex function $\eta$, there are two non-negative bounded measures $m^{\eta^{\prime\prime}}(t,x)$, $n^\eta(t,x)$ satisfying that
\begin{eqnarray}\label{2.1}
m^{\eta^{\prime\prime}}(t,x)=\int_\mR \eta^{\prime\prime}(v)m(t,x,v)dv, \ with \ m(t,x,v) \ a \ nonnegative \  measure,
\end{eqnarray}
and
\begin{eqnarray}\label{2.2}
n^\eta(t,x)=\frac{\eta^\prime(\rho(t,x))(-\Delta_x)^{\frac{\alpha}{2}}\rho(t,x)
-(-\Delta_x)^{\frac{\alpha}{2}}\eta(\rho(t,x))}{2},
\end{eqnarray}
such that the following identity holds
\begin{eqnarray}\label{2.3}
\frac{\partial}{\partial t}\eta(\rho)+\div_xQ(\rho)
+(-\Delta_x)^{\frac{\alpha}{2}}\eta(\rho)=-2m^{\eta^{\prime\prime}}(t,x)-2n^\eta(t,x),
\end{eqnarray}
in $\cD^\prime([0,T)\times\mR^d)$ with $\eta(\rho(t=0))=\eta(\rho_0)$, where $Q(\rho)=\int^\rho \eta^\prime(v)f(v)dv$.
\end{definition}

\begin{remark} (i) We define an entropy solution by the identity (\ref{2.3}), and the source or motivation for this definition comes from the $\varepsilon\rightarrow 0$ limit of the following equation directly,
\begin{eqnarray}\label{2.4}
\frac{\partial}{\partial t}\rho_\varepsilon(t,x)+\div_xF(\rho_\varepsilon)
+(-\Delta_x)^{\frac{\alpha}{2}} \rho_\varepsilon-\varepsilon\Delta\rho_\varepsilon=0.
\end{eqnarray}
  \vskip1mm\par
Indeed, if one multiplies equation (\ref{2.4}) by $ \eta^\prime(\rho_\varepsilon)$, it yields
\begin{eqnarray}\label{2.5}
\frac{\partial}{\partial t}\eta(\rho_\varepsilon)+\div_xQ(\rho_\varepsilon)
+\eta^\prime(\rho_\varepsilon)(-\Delta_x)^{\frac{\alpha}{2}}\rho_\varepsilon
=\varepsilon\eta^\prime(\rho_\varepsilon)\Delta\rho_\varepsilon.
\end{eqnarray}
With the help of the chain rule,
\begin{eqnarray}\label{2.6}
\varepsilon\eta^\prime(\rho_\varepsilon)\Delta\rho_\varepsilon=\varepsilon
\Delta\eta(\rho_\varepsilon)-\varepsilon\eta^{\prime\prime}(\rho_\varepsilon)
|\nabla\rho_\varepsilon|^2
=:\varepsilon
\Delta\eta(\rho_\varepsilon)-2m^{\eta^{\prime\prime}}_\varepsilon.
\end{eqnarray}
 \vskip1mm\par
Moreover, since $\eta$ is convex, by (\ref{1.4}),
\begin{eqnarray}\label{2.7}
\eta^\prime(\rho_\varepsilon)(-\Delta_x)^{\frac{\alpha}{2}}\rho_\varepsilon(t,x)
\geq
c_0\int_{\mR^d}\frac{\eta(\rho_\varepsilon(t,x))-
\eta(\rho_\varepsilon(t,z+x))}
{|z|^{d+\alpha}}dz
=(-\Delta_x)^{\frac{\alpha}{2}}\eta(\rho_\varepsilon(t,x)).
\end{eqnarray}
 \vskip1mm\par
Combining (\ref{2.6}) and (\ref{2.7}), we conclude from (\ref{2.5}) that
\begin{eqnarray}\label{2.8}
\frac{\partial}{\partial t}\eta(\rho_\varepsilon)+\div_xQ(\rho_\varepsilon)
+(-\Delta_x)^{\frac{\alpha}{2}}\eta(\rho_\varepsilon(t,x))
=\varepsilon
\Delta\eta(\rho_\varepsilon)-2
m^{\eta^{\prime\prime}}_\varepsilon-2n^{\eta}_\varepsilon,
\end{eqnarray}
with non-negative measures $m^{\eta^{\prime\prime}}_\varepsilon$ and $n^{\eta}_\varepsilon$. So the vanishing viscosity limit in the proceeding identity  leads to  (\ref{2.3}).
 \vskip1mm\par
(ii) Another motivation  to define entropy solutions is from \cite{CJ} Definition 2.2 and Lemma 2.4. Since
\begin{eqnarray}\label{2.9}
\|(-\Delta_x)^{\frac{\alpha}{2}}\rho(t,x)\|_{L^1(\mR^d)}
\leq C\|\rho(t)\|^{1-\alpha}_{L^1(\mR^d)}\|\rho(t)\|^\alpha_{BV(\mR^d)}
\end{eqnarray}
and
\begin{eqnarray}\label{2.10}
\|(-\Delta_x)^{\frac{\alpha}{2}}\eta(\rho(t,x))\|_{L^1(\mR^d)}
\leq C\|\rho(t)\|^{1-\alpha}_{L^1(\mR^d)}\|\rho(t)\|^\alpha_{BV(\mR^d)},
\end{eqnarray}
we have
$$
(-\Delta_x)^{\frac{\alpha}{2}}\eta(\rho)+2n^\eta(t,x)=\eta^\prime(\rho(t,x))(-\Delta_x)^{\frac{\alpha}{2}}\rho(t,x).
$$
The present definition is the same as   Definition 2.2 in \cite{CJ}. The only difference is that, here we define entropy solutions by an identity but not an inequality. As mentioned in introduction, the intermediate solution may fail to be unique, so we give an explicit formula for dissipation measure $n$, and it comes from \cite{CP} Definition 2.1 for non-isotropic degenerate parabolic-hyperbolic equation:
\begin{eqnarray}\label{2.11}
\frac{\partial}{\partial
t} \rho(t,x)+\div_xF(\rho)+\nabla\cdot(A(\rho)\nabla\rho)
=0, \ A(\rho)=\sigma(\rho)\sigma(\rho)^\top, \ \sigma\in \mR^{d\times J}.
\end{eqnarray}
 \vskip1mm\par
(iii) The main ingredient in Definition 2.1 of \cite{CP} is the chain rule for $\rho$, i.e.
\begin{eqnarray}\label{2.12}
\psi(\rho)\sigma_{i,j}(\rho)\partial_{x_i}\rho=
\partial_{x_i}\beta_{i,j}^\psi(\rho), \ \forall \ \psi\in \cD_+(\mR),
\end{eqnarray}
where $\beta_{i,j}^\psi$ is a special function, see \cite{CP}.

Even though $\rho\in L^1\cap L^\infty$ does not make the left hand side meaningful, the chain rule ensures that all manipulations legitimate in (\ref{2.12}). When the degenerate parabolic operator is replaced by a fractional operator, this chain rule may no longer hold. However, if $\rho\in L^\infty([0,T);BV(\mR^d))$ and $\alpha\in (0,1)$, with the help of (\ref{2.9})-(\ref{2.10}), for any convex smooth function $\eta$, we have
\begin{eqnarray}\label{2.13}
0\leq \frac{\eta^\prime(\rho(t,x))(-\Delta_x)^{\frac{\alpha}{2}}\rho(t,x)
-(-\Delta_x)^{\frac{\alpha}{2}}\eta(\rho(t,x))}{2}
\in L^1([0,T)\times\mR^d).
\end{eqnarray}
Note that (\ref{2.9}) and (\ref{2.10}) is meaningful if and only if $\rho\in \cC([0,T];L^1(\mR^d))\cap L^\infty([0,T);BV(\mR^d))$. If $\eta^\prime$ is bounded, then the microscopic equation (\ref{1.11}) is legitimate for  $\rho\in \cC([0,T];L^1(\mR^d))\cap L^\infty([0,T);BV(\mR^d))$. By this observation, we introduce the following definition.
\end{remark}

\begin{definition}   (Kinetic solution)  Let (\ref{1.3}) hold. A function $\rho\in \cC([0,T];L^1(\mR^d))\cap L^\infty([0,T);BV(\mR^d))$ is called a kinetic solution of (\ref{1.1})- (\ref{1.2}), if $u$, defined by (\ref{1.7}), satisfies  (\ref{1.11}), (\ref{1.9}) in $\cD^\prime([0,T)\times\mR^{d+1})$ and
 \vskip1mm\par
(i) the non-negative measure $n(t,x,v)$ is given by
\begin{eqnarray}\label{2.14}
n(t,x,v)=\frac{\sgn(\rho(t,x)-v)(-\Delta_x)^{\frac{\alpha}{2}}\rho(t,x)
-(-\Delta_x)^{\frac{\alpha}{2}}|\rho(t,x)-v|}{2};
\end{eqnarray}
\vskip1mm\par
(ii) the nonnegative measure $m+n$  fulfils the condition
\begin{eqnarray}\label{2.15}
\int_0^T\int_{\mR^d}(m+n)(dt,dx,v)\in L^\infty_0(\mR).
\end{eqnarray}
\end{definition}

\begin{remark} Note that $\rho\in \cC([0,T];L^1(\mR^d))$, $n(t,x,v)\in L^\infty_v(\mR;L^1_{t,x}([0,T)\times\mR^d))$. Hence the nonnegative measure $m$ in Definition 2.2 is continuous in $t$ in the sense  that
$$
\lim_{s\rightarrow t} \int_0^s\varphi(x)\phi(v)m(ds,dx,dv)=\int_0^t\varphi(x)\phi(v)m(ds,dx,dv),
$$
for   $\varphi\in \cD(\mR^d)$ and $\phi\in\cD(\mR)$, which imply that the preceding  definition is equivalent to
\begin{eqnarray*}
&&\int_{\mR^{d+1}}u(t,x,v)\varphi(x)\phi(v)dxdv-\int_{\mR^{d+1}}\chi_{\rho_0}(v)\varphi(x)\phi(v)dxdv
\cr\cr&=&\int^t_0\int_{\mR^{d+1}}u(s,x,v)f(v)\cdot\nabla\varphi(x)\psi(v)dxdvds
-\int^t_0\int_{\mR^{d+1}}u(s,x,v)\cL\varphi(x)\psi(v)dxdvds\cr\cr&&-
\int^t_0\int_{\mR^{d+1}}\varphi(x)\frac{\partial}{\partial v}\psi(v)(m+n)(ds,dx,dv),
\end{eqnarray*}
for   $t\in (0,T)$,  $\varphi\in \cD(\mR^d)$  and  $\phi\in\cD(\mR)$.
\end{remark}
 \vskip1mm\par
Now, we are in a position to show the relationship between entropy
solutions and kinetic solutions for (\ref{1.1})- (\ref{1.2}).
\begin{theorem} \textbf{(Kinetic formulation)}  Let (\ref{1.3}) be valid, $\rho_0\in L^\infty(\mR^d)$ and $u(t,x,v)=\chi_{\rho}(v)$.
  \vskip1mm\par
(i) If $\rho$ is an entropy solution of (\ref{1.1})-(\ref{1.2}), then
it is also a kinetic solution. Besides, the nonnegative measures $m$ and $n$ are bounded and supported in $[0,T]\times \mR^d\times [-M,M]$ for $M=\|\rho\|_{L^\infty([0,T)\times\mR^d)}$, and further satisfy  (\ref{2.15}).
 \vskip1mm\par
(ii) If $\rho$ is a kinetic
solution of (\ref{1.1})-(\ref{1.2}), then it is an entropy solution as well.
\end{theorem}
\textbf{Proof.} By the following relationship:
$$
\int_{\mR}S^\prime(v)u(t,x,v)dv=S(\rho(t,x)), \ \forall \ S\in\cC^1(\mR),
$$
we clearly get the conclusion (ii).  It remains  to verify (i).
 \vskip1mm\par
Indeed, if $\rho$ is an entropy solution, then from (\ref{2.3}), by an approximation, we deduce that
\begin{eqnarray}\label{2.16}
\frac{\partial }{\partial t}\eta(\rho,v) +\div_xQ(\rho,v)
+(-\Delta_x)^{\frac{\alpha}{2}}|\rho-v|
=-2m-2n,
\end{eqnarray}
where $\eta(\rho,v)=|\rho-v|, Q(\rho,v)=\sgn(\rho-v)[F(\rho)-F(v)]$.
 \vskip1mm\par
By differentiating (\ref{2.16}) in $v$ in the distributions sense, we obtain the equation (\ref{1.11}). Besides, from (\ref{2.16}), $m$ and $n$ are nonnegative and supported in $[0,T)\times \mR^d\times [-M,M]$.
 \vskip1mm\par
Furthermore, if one integrates the identity (\ref{2.16}) in $x$ on $\mR^d$, then
\begin{eqnarray*}
\int_0^T\int_{\mR^d}(m+n)(dt,dx,v)
=
\frac{1}{2}\int_{\mR^d}[|\rho_0(x)-v|-|\rho(T,x)-v|]dx
\leq\sup_{0\leq t\leq T}  \int_{\mR^d}|\rho(t,x)|dx.
\end{eqnarray*}
Therefore $m+n$ is bounded and (\ref{2.15}) holds. We complete the proof.

\begin{remark}  The proof here is   analogue to  that   for
\begin{eqnarray*}
\frac{\partial}{\partial
t} \rho(t,x)+\div_xF(\rho)=0, \ \ (t,x)\in (0,T) \times \mR^d,
\end{eqnarray*}
in \cite{LPT}; so we omit some   details.
 \vskip1mm\par
(ii) Observe that
$$
(\rho-v)_+=\frac{|\rho-v|+(\rho-v)}{2}, \ \ (\rho-v)_-=\frac{|\rho-v|-(\rho-v)}{2},
$$
Thus  if $\rho$ is an entropy solution,  we  can take entropy-entropy flux pairs by
$(\eta(\rho,v)=(\rho-v)_+, Q(\rho,v)=\sgn(\rho-v)_+[F(\rho)-F(v)])$ and $(\eta(\rho,v)=(\rho-v)_-, Q(\rho,v)=\sgn(\rho-v)_-[F(\rho)-F(v)])$ respectively,  and we can  estimate  that
\begin{eqnarray}\label{2.17}
\int_0^T\int_{\mR^d}(m+n)(dt,dx,v)\leq
\frac{1}{2}\int_{\mR^d}
|(\rho_0(x)-v)_+-(\rho(T,x)-v)_+|
dx
\end{eqnarray}
and
\begin{eqnarray}\label{2.18}
\int_0^T\int_{\mR^d}(m+n)(dt,dx,v)\leq
\frac{1}{2}\int_{\mR^d}|(\rho_0(x)-v)_--(\rho(T,x)-v)_-|dx.
\end{eqnarray}
From (\ref{2.17})-(\ref{2.18}), we derive
\begin{eqnarray*}
\lim_{v\rightarrow \pm\infty}\int_0^T\int_{\mR^d}(m+n)(dt,dx,v)=0,
\end{eqnarray*}
i.e. (\ref{2.15}) is true. Observing that the right hand sides in (\ref{2.17}) and (\ref{2.18}) are meaningful   if $\rho$ is a kinetic solution. Thus in Definition 2.2, we add the condition (\ref{2.15}).
\end{remark}

\begin{remark} The preceding result holds  as well for the non-homogeneous fractional convection-diffusion problem
\begin{eqnarray}\label{2.19}
\left\{ \begin{array}{ll} \frac{\partial}{\partial
t} \rho(t,x)+\div_xF(\rho)+(-\Delta_x)^{\frac{\alpha}{2}}
B(\rho)=A(\rho), \ \ (t,x)\in (0,T) \times \mR^d,
\\ \rho(t=0,x)=\rho_0(x), \ \ x\in \mR^d,
\end{array}\right.
\end{eqnarray}
if
\begin{eqnarray}\label{2.20}
A(0)=0, \ A\in W^{1,1}_{loc}(\mR), \ B(0)=0, \ B\in W^{1,\infty}_{loc}(\mR), \  B^\prime\geq 0.
\end{eqnarray}
But now the Cauchy problem (\ref{1.9}) with (\ref{1.11}) should be replaced  by
\begin{eqnarray}\label{2.21}
\left\{ \begin{array}{ll} \frac{\partial}{\partial t}u+f(v)\cdot \nabla_xu+ A(v)\frac{\partial}{\partial v}u+b(v)(-\Delta_x)^{\frac{\alpha}{2}} u=\frac{\partial}{\partial v}(m+n), \ (t,x,v)\in (0,T) \times \mR^d\times \mR, \\ u(t=0)=\chi_{\rho_0}(v), \  (x,v)\in \mR^d\times \mR, \end{array}\right.
\end{eqnarray}
where $b(v)=B^\prime(v)$,
\begin{eqnarray}\label{2.22}
n(t,x,v)=\frac{\sgn(\rho(t,x)-v)(-\Delta_x)^{\frac{\alpha}{2}}B(\rho(t,x))
-(-\Delta_x)^{\frac{\alpha}{2}}|B(\rho(t,x))-B(v)|}{2}.
\end{eqnarray}
\end{remark}

\section{Uniqueness and existence of kinetic solutions}\label{sec3}\setcounter{equation}{0}

In this section, we are interested in the Cauchy problem (\ref{1.1})-(\ref{1.2}) and it is ready for us to state our main result.
\begin{theorem}
Let (\ref{1.3}) hold. Then there is a unique kinetic solution of the nonlocal Cauchy problem (\ref{1.1})-(\ref{1.2}).
\end{theorem}
     \vskip0mm\noindent
\textbf{Proof.} \textbf{(Uniqueness)} Let $\rho_i \ (i=1,2)$ be kinetic solutions of (\ref{1.1})-(\ref{1.2}). Then  for both $u_i=\chi_{\rho_i} \ (i=1,2)$
\begin{eqnarray}\label{3.1}
\left\{ \begin{array}{ll}\frac{\partial }{\partial t}u_i(t,x,v)+f(v)\cdot \nabla_xu_i+(-\Delta_x)^{\frac{\alpha}{2}}u_i=\frac{\partial
}{\partial v}m_i(t,x,v)+\frac{\partial
}{\partial v}n_i(t,x,v),  \\ u_i(t=0)=\chi_{\rho_0}(v) \in L^1(\mR^{d+1}),
\end{array}\right.
\end{eqnarray}
with the nonnegative measures $m_i,
n_i$ satisfying (\ref{2.14}) and (\ref{2.15}).
 \vskip1mm\par
We set
\begin{eqnarray}\label{3.2}
n_i^1(t,x,v)=\sgn(\rho_i(t,x)-v)(-\Delta_x)^{\frac{\alpha}{2}}\rho_i(t,x),
\ \ n_i^2(t,x,v)=-(-\Delta_x)^{\frac{\alpha}{2}}|\rho_i(t,x)-v|.
\end{eqnarray}
Then
\begin{eqnarray}\label{3.3}
n_i^1,  n_i^2\in L^\infty_0(\mR_v;L^1([0,T)\times\mR_x^d)) \ \mbox{and} \ 2n_i=n_i^1+n_i^2.
\end{eqnarray}
 \vskip1mm\par
For $\varepsilon_1, \varepsilon_2, \sigma>0,$   define
$$
\varrho_{1,\varepsilon_1}(t)=\frac{1}{\varepsilon_1}
\varrho_1(\frac{t}{\varepsilon_1}), \ \varrho_{2,\varepsilon_2}(x)=\frac{1}{\varepsilon_2^d}
\varrho_2(\frac{x}{\varepsilon_2}), \ \varrho_{3,\sigma}(v)=\frac{1}{\sigma}
\varrho_3(\frac{v}{\sigma}),
$$
here $\varrho_1$, $\varrho_2$ and $\varrho_3$ are three nonnegative normalized regularizing kernels, satisfying
$$
\mbox{supp} \varrho_1\subset (-1,0), \ \mbox{supp}\varrho_2\subset B_1(0), \ \mbox{supp} \varrho_3\subset (-1,1).
$$
Then $u_{i,\varepsilon}:=u_i\ast\varrho_{1,\varepsilon_1}\ast
\varrho_{2,\varepsilon_2} \  (i=1,2$) yield
\begin{eqnarray}\label{3.4}
\left\{ \begin{array}{ll} \frac{\partial }{\partial t}u_{i,\varepsilon}+f(v)\cdot \nabla_xu_{i,\varepsilon}+ (-\Delta_x)^{\frac{\alpha}{2}}u_{i,\varepsilon}
=\frac{\partial
}{\partial v}m_{i,\varepsilon}(t,x,v)+\frac{\partial
}{\partial v}n_{i,\varepsilon}(t,x,v), \\  u_{i,\varepsilon}(t=0)=
\chi_{\rho_0}\ast\varrho_{2,\varepsilon_2}(x),
\end{array}\right.
\end{eqnarray}
and $u^\sigma_{i,\varepsilon}:=u_i\ast\varrho_{1,\varepsilon_1}\ast
\varrho_{2,\varepsilon_2}\ast\varrho_{3,\sigma} \  (i=1,2$) fulfill
\begin{eqnarray}\label{3.5}
\left\{ \begin{array}{ll}\frac{\partial }{\partial t}u^\sigma_{i,\varepsilon}+f(v)\cdot \nabla_xu^\sigma_{i,\varepsilon}+(-\Delta_x)^{\frac{\alpha}{2}}u^\sigma_{i,\varepsilon}
=\frac{\partial
}{\partial v}m_{i,\varepsilon}^\sigma+\frac{\partial
}{\partial v}n_{i,\varepsilon}^\sigma+R_{i,\varepsilon}^{\sigma}, \\ u^\sigma_{i,\varepsilon}(t=0)=
\chi_{\rho_0}\ast\varrho_{2,\varepsilon_2}(x)\ast\varrho_{3,\sigma}(v),
\end{array}\right.
\end{eqnarray}
with
\begin{eqnarray}\label{3.6}
 R_{i,\varepsilon}^{\sigma}=f(v)\cdot \nabla_xu^\sigma_{i,\varepsilon}-(f(v)\cdot \nabla_xu_{i,\varepsilon})\ast\varrho_{3,\sigma},
\end{eqnarray}
here we define $u_i(t,x,v):=0$, when $t\bar{\in}[0,T)$.
\vskip1mm\par
In view of $|u_{i,\varepsilon}|=\sgn(v)u_{i,\varepsilon}$, we get from (\ref{3.4})-(\ref{3.6}) that
\begin{eqnarray*}
\frac{\partial }{\partial t}|u_{i,\varepsilon}|+f(v)\cdot \nabla_x|u_{i,\varepsilon}|+(-\Delta_x)^{\frac{\alpha}{2}}|u_{i,\varepsilon}|= \sgn(v)[\frac{\partial
}{\partial v}m_{i,\varepsilon}+\frac{\partial
}{\partial v}n_{i,\varepsilon}],
\end{eqnarray*}
and
\begin{eqnarray*}
&&\frac{\partial}{\partial t}(u^\sigma_{1,\varepsilon}u^\sigma_{2,\varepsilon}) +f\cdot \nabla_x(u^\sigma_{1,\varepsilon}u^\sigma_{2,\varepsilon})+[u^\sigma_{2,\varepsilon}
(-\Delta_x)^{\frac{\alpha}{2}}u^\sigma_{1,\varepsilon}+
u^\sigma_{1,\varepsilon}
(-\Delta_x)^{\frac{\alpha}{2}}u^\sigma_{2,\varepsilon}]
\cr\cr&=&u^\sigma_{1,\varepsilon}\frac{\partial
}{\partial v}[m^\sigma_{2,\varepsilon}+n^\sigma_{2,\varepsilon}]+u^\sigma_{2,\varepsilon}
\frac{\partial
}{\partial v}[m^\sigma_{1,\varepsilon}+n^\sigma_{1,\varepsilon}]+
u^\sigma_{1,\varepsilon}R_{2,\varepsilon}^{\sigma}
+u^\sigma_{2,\varepsilon}R_{1,\varepsilon}^{\sigma}.
\end{eqnarray*}
Hence
\begin{eqnarray}\label{3.7}
&&\!\!\!\frac{\partial }{\partial t}[|u_{1,\varepsilon}|\!+\!|u_{2,\varepsilon}|\!-\!2u^\sigma_{1,\varepsilon}u^\sigma_{2,\varepsilon}]+f(v)\cdot \nabla_x[|u_{1,\varepsilon}|\!+\!|u_{2,\varepsilon}|\!-\!2u^\sigma_{1,\varepsilon}u^\sigma_{2,\varepsilon}]
+(-\Delta_x)^{\frac{\alpha}{2}}[|u_{1,\varepsilon}|
+|u_{2,\varepsilon}|]
\cr\cr&=&\!\!\!I_1+I_2-2u^\sigma_{1,\varepsilon}R_{2,\varepsilon}^{\sigma}
-2u^\sigma_{2,\varepsilon}R_{1,\varepsilon}^{\sigma},
\end{eqnarray}
where
\begin{eqnarray}\label{3.8}
I_1(t,x,v)=\sgn(v)\frac{\partial
}{\partial v}[m_{1,\varepsilon}
+m_{2,\varepsilon}]-2[u^\sigma_{1,\varepsilon}\frac{\partial
}{\partial v}m_{2,\varepsilon}^\sigma+u^\sigma_{2,\varepsilon}\frac{\partial
}{\partial v}m_{1,\varepsilon}^\sigma]
\end{eqnarray}
and
\begin{eqnarray}\label{3.9}
I_2(t,x,v)=\sgn(v)\frac{\partial
}{\partial v}[n_{1,\varepsilon}
+n_{2,\varepsilon}]
-2[u^\sigma_{1,\varepsilon}\frac{\partial
}{\partial v}n_{2,\varepsilon}^{1,\sigma}+u^{\sigma}_{2,\varepsilon}\frac{\partial
}{\partial v}n^{1,\sigma}_{1,\varepsilon}].
\end{eqnarray}
  \vskip1mm\par
Let $\theta$ and $\xi$ be two cut-off functions, with variables $x$ and $v$ respectively, i.e. $\theta\in{\cD}(\mR^d)$, $\xi\in{\cD}(\mR)$,
\begin{eqnarray}\label{3.10}
0\leq \theta, \xi \leq 1, \ \ \theta=\cases {1,  \ \   |x|\leq 1,   \cr 0, \ \  |x|\geq 2,}  \ \
\xi=\cases {1,  \ \   |v|\leq 1,   \cr 0, \ \  |v|\geq 2,} \ \ v\xi^\prime(v)\leq 0,
\end{eqnarray}
and for $p,k\in\mN$, we denote by  $\theta_p(x)=\theta(\frac{x}{p})$ and $\xi_k(v)=\xi(\frac{v}{k})$.
 \vskip1mm\par
Now let us estimate the right hand sides in (\ref{3.7}). Initially, we have the following estimate for the last two error terms,
\begin{eqnarray}\label{3.11}
\lim_{\sigma\rightarrow 0}\int_{\mR^{d+1}}[u^\sigma_{1,\varepsilon}R_{2,\varepsilon}^{\sigma}
+u^\sigma_{2,\varepsilon}R_{1,\varepsilon}^{\sigma}]\xi_k(v)\theta_p(x)dxdv=0
\end{eqnarray}
for fixed $\varepsilon_1,\varepsilon_2, k$ and $p$. It remains to reckon the others.
 \vskip1mm\par
Note that
\begin{eqnarray}\label{3.12}
&&\lim_{\sigma\rightarrow 0}\int_{\mR^{d+1}}\partial_vu^\sigma_{1,\varepsilon}(t,x,v)m^\sigma_{2,\varepsilon}(t,x,v)
\xi_k(v)\theta_p(x)dxdv\cr\cr&=&\int_{\mR^{d+1}}\partial_v
u_{1,\varepsilon}(t,x,v)m_{2,\varepsilon}(t,x,v)
\xi_k(v)\theta_p(x)dxdv\cr\cr&=&\int_{\mR^{d+1}}\partial_v
\int_0^t\int_{\mR^d}\chi_{\rho_1(s,y)}(v)\varrho_{1,\varepsilon_1}(t-s)
\varrho_{2,\varepsilon_2}(x-y)dydsm_{2,\varepsilon}(t,x,v)
\xi_k(v)\theta_p(x)dxdv
\cr\cr&=&\int_{\mR^{d+1}}
\int_0^t\int_{\mR^d}[\delta(v)-\delta(v-\rho_1(s,y))]\varrho_{1,\varepsilon_1}(t-s)
\varrho_{2,\varepsilon_2}(x-y)dydsm_{2,\varepsilon}(t,x,v)
\xi_k(v)\theta_p(x)dxdv
\cr\cr
&\leq& \int_{\mR^d}m_{2,\varepsilon}(t,x,0)\theta_p(x)dx,
\end{eqnarray}
where $\delta$ is the Dirac mass concentrated at 0.
 \vskip1mm\par
Similar calculations also lead to
\begin{eqnarray}\label{3.13}
\lim_{\sigma\rightarrow 0}\int_{\mR^{d+1}}\partial_vu^\sigma_{2,\varepsilon}(t,x,v)m^\sigma_{1,\varepsilon}(t,x,v)
\xi_k(v)\theta_p(x)dxdv
\leq \int_{\mR^d}m_{1,\varepsilon}(t,x,0)\theta_p(x)dx.
\end{eqnarray}
Therefore
\begin{eqnarray}\label{3.14}
\!\!&&\lim_{\sigma\rightarrow 0}\int_{\mR^{d+1}}I_1(t,x,v)\xi_k(v)\theta_p(x)dxdv \cr\cr\!\!&\leq&\!\!
2\int_{\mR^{d+1}}
[u_{1,\varepsilon}m_{2,\varepsilon}+u_{2,\varepsilon}m_{1,\varepsilon}]
\partial_v\xi_k\theta_pdxdv
\!-\!2\int_{\mR^{d+1}}[m_{1,\varepsilon}+m_{2,\varepsilon}]\theta_p
\sgn(v) \partial_v\xi_k dxdv.
\end{eqnarray}
  \vskip1mm\par
From (\ref{3.10}), one can deduce that $\sgn(v)\xi^\prime(v)\leq 0$.
By virtue of (\ref{1.7}) and (\ref{3.10}), it follows that
$$
u_{1,\varepsilon}(t,x,v)\partial_v\xi_k(v)=
\frac{1}{k}\int_0^t\int_{\mR^d}\chi_{\rho_1(s,y)}(v)\xi^\prime(\frac{v}{k})\varrho_{1,\varepsilon_1}(t-s)
\varrho_{2,\varepsilon_2}(x-y)dyds\leq0
$$
and
$$
u_{2,\varepsilon}(t,x,v)\partial_v\xi_k(v)=
\frac{1}{k}\int_0^t\int_{\mR^d}\chi_{\rho_2(s,y)}(v)\xi^\prime(\frac{v}{k})\varrho_{1,\varepsilon_1}(t-s)
\varrho_{2,\varepsilon_2}(x-y)dyds\leq0.
$$
Moreover, since $m_{1,\varepsilon}, m_{2,\varepsilon}\geq 0$, from (\ref{3.14}), it leads to
\begin{eqnarray*}
\lim_{\sigma\rightarrow 0}\int_{\mR^{d+1}}I_1(t,x,v)\xi_k(v)\theta_p(x)dxdv \leq
-2\int_{\mR^{d+1}}[m_{1,\varepsilon}+m_{2,\varepsilon}]\theta_p(x)
\sgn(v) \partial_v\xi_k(v) dxdv.
\end{eqnarray*}
In view of (\ref{2.15}),  we  have
\begin{eqnarray}\label{3.15}
\lim_{k\rightarrow \infty}\lim_{p\rightarrow \infty}\lim_{\sigma\rightarrow 0}\int_{\mR^{d+1}}I_1(t,x,v)\xi_k(v)\theta_p(x)dxdv \leq 0.
\end{eqnarray}
  \vskip1mm\par
Now let us estimate the term $I_2$ and firstly, via   integration by parts,
\begin{eqnarray}\label{3.16}
&&-2\int_{\mR^{d+1}}u_{1,\varepsilon}^\sigma(t,x,v)\partial_v
n_{2,\varepsilon}^{1,\sigma}(t,x,v)
\xi_k(v)\theta_p(x)dxdv\cr\cr
&=&2\int_{\mR^{d+1}}\partial_vu_{1,\varepsilon}^\sigma
n_{2,\varepsilon}^{1,\sigma}
\theta_p(x)\xi_k(v)dxdv+
2\int_{\mR^{d+1}}u_{1,\varepsilon}^\sigma(t,x,v)n_{2,\varepsilon}^{1,\sigma}
(t,x,v)
\partial_v\xi_k(v)\theta_p(x)dxdv
\cr\cr&=&2
\int_{\mR^{d+1}}\partial_v
\int_0^t\int_{\mR^{d+1}}\chi_{\rho_1(s,y)}(v-z)\varrho_{1,\varepsilon_1}(t-s)
\varrho_{2,\varepsilon_2}(x-y)\varrho_{3,\sigma}(z)dydsdz
\cr\cr&&\times n_{2,\varepsilon}^{1,\sigma}(t,x,v)
\theta_p(x)\xi_k(v)dxdv+
2\int_{\mR^{d+1}}u_{1,\varepsilon}^\sigma(t,x,v)n_{2,\varepsilon}^{1,\sigma}(t,x,v)
\partial_v\xi_k(v)\theta_p(x)dxdv
\cr\cr&=&2\int_{\mR^{d+1}}
\int_0^t\int_{\mR^{d+1}}[\delta(v-z)-\delta(v-z-\rho_1(s,y))]\varrho_{1,\varepsilon_1}(t-s)
\varrho_{2,\varepsilon_2}(x-y)\varrho_{3,\sigma}(z)dydsdz
\cr\cr&&\times n_{2,\varepsilon}^{1,\sigma}
\theta_p(x)\xi_k(v)dxdv+
2\int_{\mR^{d+1}}u_{1,\varepsilon}^\sigma(t,x,v)n_{2,\varepsilon}^{1,\sigma}(t,x,v)
\partial_v\xi_k(v)\theta_p(x)dxdv
\cr\cr
&=& 2\int_{\mR^{d+1}}n_{2,\varepsilon}^{1,\sigma}(t,x,z)\varrho_{3,\sigma}(z)
\theta_p(x)dxdz\cr\cr
&&-
2\int_{\mR^{d+1}}[
n_{2,\varepsilon}^{1,\sigma}(t,x,\rho_1(\cdot,\cdot)+z)\xi_k(\rho_1(\cdot,\cdot)+z)]\varrho_{3,\sigma}(z)
\ast\varrho_{1,\varepsilon_1}\ast\varrho_{2,\varepsilon_2}(t,x)
\theta_p(x)dxdz
\cr\cr&&+
2\int_{\mR^{d+1}}u_{1,\varepsilon}^\sigma(t,x,v)n^{1,\sigma}_{2,\varepsilon}(t,x,v)
\partial_v\xi_k(v)\theta_p(x)dxdv.
\end{eqnarray}
\vskip1mm\par
An analogue calculation also implies that
\begin{eqnarray}\label{3.17}
&&-2\int_{\mR^{d+1}}u_{2,\varepsilon}^\sigma(t,x,v)\partial_v
n_{1,\varepsilon}^{1,\sigma}(t,x,v)
\xi_k(v)\theta_p(x)dxdv\cr\cr
&=& 2\int_{\mR^{d+1}}n_{1,\varepsilon}^{1,\sigma}(t,x,z)\varrho_{3,\sigma}(z)
\theta_p(x)dxdz\cr\cr
&&-
2\int_{\mR^{d+1}}[
n_{1,\varepsilon}^{1,\sigma}(t,x,\rho_2(\cdot,\cdot)+z)\xi_k(\rho_2(\cdot,\cdot)+z)]\varrho_{3,\sigma}(z)
\ast\varrho_{1,\varepsilon_1}\ast\varrho_{2,\varepsilon_2}(t,x)
\theta_p(x)dxdz
\cr\cr&&+
2\int_{\mR^{d+1}}u_{2,\varepsilon}^\sigma(t,x,v)n^{1,\sigma}_{1,\varepsilon}(t,x,v)
\partial_v\xi_k(v)\theta_p(x)dxdv.
\end{eqnarray}
  \vskip1mm\par
By (\ref{3.9}), (\ref{3.16})-(\ref{3.17}), we get
\begin{eqnarray}\label{3.18}
&&\lim_{\sigma\rightarrow 0}\int_{\mR^{d+1}}I_2(t,x,v)\xi_k(v)\theta_p(x)dxdv
\cr\cr&=&-2\int_{\mR^d}[n_{2,\varepsilon}^{2}(t,x,0)+
n_{1,\varepsilon}^{2}(t,x,0)]
\theta_pdx-2\int_{\mR^{d+1}}
[n_{1,\varepsilon}(t,x,v)+
n_{2,\varepsilon}(t,x,v)]\sgn(v)
\partial_v\xi_k\theta_pdxdv\cr\cr&&-
2\int_{\mR^d}[
n_{2,\varepsilon}^{1}(t,x,\rho_1(\cdot,\cdot))\xi_k(\rho_1(\cdot,\cdot))+
n_{1,\varepsilon}^{1}(t,x,\rho_2(\cdot,\cdot))
\xi_k(\rho_2(\cdot,\cdot))]
\ast\varrho_{1,\varepsilon_1}\ast\varrho_{2,\varepsilon_2}(t,x)
\theta_p(x)dx
\cr\cr&&+
2\int_{\mR^{d+1}}[u_{1,\varepsilon}(t,x,v)n^{1}_{2,\varepsilon}(t,x,v)+
u_{2,\varepsilon}(t,x,v)n^{1}_{1,\varepsilon}(t,x,v)]
\partial_v\xi_k(v)\theta_p(x)dxdv
\cr\cr&=&-2\int_{\mR^d}[|\rho_1(t,x)|_{\varepsilon}+
|\rho_2(t,x)|_{\varepsilon}](-\Delta_x)^{\frac{\alpha}{2}}\theta_p(x)dx
\cr\cr&&-2\int_{\mR^{d+1}}
[n_{1,\varepsilon}(t,x,v)+
n_{2,\varepsilon}(t,x,v)]\sgn(v)
\partial_v\xi_k(v)\theta_p(x)dxdv
\cr\cr&&
-2\int_{\mR^d}[
n_{2,\varepsilon}^{1}(t,x,\rho_1(\cdot,\cdot))\xi_k(\rho_1(\cdot,\cdot))+
n_{1,\varepsilon}^{1}(t,x,\rho_2(\cdot,\cdot))
\xi_k(\rho_2(\cdot,\cdot))]
\ast\varrho_{1,\varepsilon_1}\ast\varrho_{2,\varepsilon_2}(t,x)
\theta_p(x)dx
\cr\cr&&+
2\int_{\mR^{d+1}}[u_{1,\varepsilon}(t,x,v)n^{1}_{2,\varepsilon}(t,x,v)+
u_{2,\varepsilon}(t,x,v)n^{1}_{1,\varepsilon}(t,x,v)]
\partial_v\xi_k(v)\theta_p(x)dxdv.
\end{eqnarray}
For $\varepsilon_1,\varepsilon_2$   fixed,
$$
[u_{1,\varepsilon}(t,x,v)n^{1}_{2,\varepsilon}(t,x,v)+
u_{2,\varepsilon}(t,x,v)n^{1}_{1,\varepsilon}(t,x,v)]
\in \cC([0,T];L^1(\mR^{d+1})).
$$
So
\begin{eqnarray}\label{3.19}
\lim_{k\rightarrow \infty}\lim_{p\rightarrow \infty}\int_{\mR^{d+1}}[u_{1,\varepsilon}(t,x,v)n^{1}_{2,\varepsilon}(t,x,v)+
u_{2,\varepsilon}(t,x,v)n^{1}_{1,\varepsilon}(t,x,v)]
\partial_v\xi_k(v)\theta_p(x)dxdv=0.
\end{eqnarray}
By (\ref{3.3}),
\begin{eqnarray}\label{3.20}
\lim_{k\rightarrow \infty}\lim_{p\rightarrow \infty}\int_{\mR^{d+1}}
[n_{1,\varepsilon}(t,x,v)+
n_{2,\varepsilon}(t,x,v)]\sgn(v)
\partial_v\xi_k(v)\theta_p(x)dxdv=0.
\end{eqnarray}
Combining (\ref{3.19}) and (\ref{3.20}), from (\ref{3.18}), we assert that
\begin{eqnarray}\label{3.21}
&&\lim_{k\rightarrow \infty}\lim_{p\rightarrow \infty}\lim_{\sigma\rightarrow 0}\int_{\mR^{d+1}}I_2(t,x,v)\xi_k(v)\theta_p(x)dxdv
\cr\cr&=&-2\int_{\mR^d}[
n_{2,\varepsilon}^{1}(t,x,\rho_1(\cdot,\cdot))+
n_{1,\varepsilon}^{1}(t,x,\rho_2(\cdot,\cdot))]
\ast\varrho_{1,\varepsilon_1}\ast\varrho_{2,\varepsilon_2}(t,x)
dx.
\end{eqnarray}
 \vskip1mm\par
From (\ref{3.21}), if we take $\varepsilon_1\downarrow0$, $\varepsilon_2\downarrow0$ in turn, then
\begin{eqnarray}\label{3.22}
&&\lim_{\varepsilon_2\rightarrow 0}\lim_{\varepsilon_1\rightarrow 0}\lim_{k\rightarrow \infty}\lim_{p\rightarrow \infty}\lim_{\sigma\rightarrow 0}\int_{\mR^{d+1}}I_2(t,x,v)\xi_k(v)\theta_p(x)dxdv
\cr\cr&=&-2\int_{\mR^d}\sgn(\rho_2(t,x)-\rho_1(t,x))[
(-\Delta_x)^{\frac{\alpha}{2}}\rho_2(t,x)-(-\Delta_x)^{\frac{\alpha}{2}}\rho_1(t,x)
]dx
\cr\cr&\leq&
-2\int_{\mR^d}(-\Delta_x)^{\frac{\alpha}{2}}|\rho_2(t,x)-\rho_1(t,x)|
dx\cr\cr&=&0.
\end{eqnarray}
 \vskip1mm\par
By Remark 2.3, then (\ref{3.7}) implies
\begin{eqnarray}\label{3.23}
&& \frac{d}{dt}\int_{\mR^{d+1}}[|u_{1,\varepsilon}|+|u_{2,\varepsilon}|-
2u^\sigma_{1,\varepsilon}u^\sigma_{2,\varepsilon}]
\xi_k(v)\theta_p(x)dxdv
\cr\cr&=&\int_{\mR^{d+1}}[|u_{1,\varepsilon}|+|u_{2,\varepsilon}|-
2u^\sigma_{1,\varepsilon}u^\sigma_{2,\varepsilon}]
\xi_k(v)f(v)\cdot\nabla_x\theta_p(x)dxdv
\cr\cr&&
-\int_{\mR^{d+1}}(-\Delta_x)^{\frac{\alpha}{2}}\theta_p(x)
[|u_{1,\varepsilon}|+|u_{2,\varepsilon}|]\xi_k(v)dxdv
-
2\int_{\mR^{d+1}}[u^\sigma_{1,\varepsilon}R_{2,\varepsilon}^{\sigma}
+u^\sigma_{2,\varepsilon}R_{1,\varepsilon}^{\sigma}]\xi_k\theta_pdxdv
\cr\cr&&+\int_{\mR^{d+1}}I_1(t,x,v)\xi_k(v)\theta_p(x)dxdv
+\int_{\mR^{d+1}}I_2(t,x,v)\xi_k(v)\theta_p(x)dxdv.
\end{eqnarray}
According to (\ref{3.11}), (\ref{3.15}) and (\ref{3.22}), if we let $\sigma\downarrow0$ first, $p\uparrow\infty$ second, $k\uparrow\infty$ third, $\varepsilon_1\downarrow0$  fourth and $\varepsilon_2\downarrow0$ last,  we  conclude  from (\ref{3.23}) that
\begin{eqnarray}\label{3.24}
\frac{d}{dt}\int_{\mR^{d+1}}[|u_{1}|+|u_{2}|-2u_{1}u_{2}]
dxdv\leq 0.
\end{eqnarray}
Since $|u_1|=|u_1|^2, |u_2|=|u_2|^2, |u_1-u_2|=|u_1-u_2|^2$, from (\ref{3.24}), we end up with
\begin{eqnarray*}
\frac{d}{dt} \int_{\mR^{d+1}}|u_1(t)-u_2(t)|dxdv \leq 0,
\end{eqnarray*}
which indicates
\begin{eqnarray}\label{3.25}
\int_{\mR^d}|\rho_1(t)-\rho_2(t)|dx=\int_{\mR^{d+1}}|u_1(t)-u_2(t)|
dxdv=0.
\end{eqnarray}
From this, we finish the proof for the uniqueness.
 \vskip2mm\noindent
\textbf{(Existence)} We prove the existence by a vanishing viscosity method.  Assume that  $\rho_0\in L^\infty\cap L^1\cap BV(\mR^d)$.
 \vskip1mm\par
Consider the Cauchy problem:
\begin{eqnarray}\label{3.26}
\left\{ \begin{array}{ll}
\frac{\partial }{\partial t}\rho_\varepsilon(t,x)+\div_xF(\rho_\varepsilon)
+(-\Delta_x)^{\frac{\alpha}{2}}\rho_\varepsilon-
\varepsilon\Delta\rho_\varepsilon=0, \ \ (t,x)\in (0,T) \times \mR^d,
\\ \rho_\varepsilon(t=0)=\rho_0, \ \ x\in\mR^d.
\end{array}\right.
\end{eqnarray}
With the classical parabolic theory (see \cite{VH}), there is a unique strong solution $\rho_\varepsilon$ of (\ref{3.26}) and for any smooth convex function $\eta$, (\ref{2.8}) holds (see Remark 2.1). Moreover, the following inequalities hold
\begin{eqnarray}\label{3.27}
\left\{ \begin{array}{ll} \|\rho_\varepsilon(t)\|_{L^1(\mR^d)}\leq  \|\rho_0\|_{L^1(\mR^d)}, \  \|\rho_\varepsilon(t)\|_{BV(\mR^d)}\leq  \|\rho_0\|_{BV(\mR^d)},
 \\
\|\partial_t\rho_\varepsilon(t)\|_{L^1(\mR^d)}\leq  C(\|\rho_0\|_{BV(\mR^d)}+\|\rho_0\|_{L^1(\mR^d)}).
\end{array}\right.
\end{eqnarray}
   \vskip1mm\par
Indeed, if we choose $\eta(\rho)=|\rho|$, with the help of entropy inequality (\ref{2.3}) (since a classical solution is also an entropy solution), it follows that
\begin{eqnarray}\label{3.28}
\frac{\partial}{\partial t}|\rho_\varepsilon|+\div_x
[\sgn(\rho_\varepsilon)F(\rho_\varepsilon)]
\leq
\varepsilon\Delta|\rho_\varepsilon|-(-\Delta_x)^{\frac{\alpha}{2}}
|\rho_\varepsilon|.
\end{eqnarray}
  \vskip1mm\par
Integrating both hand sides of (\ref{3.28}) on $\mR^d$, we  obtain
\begin{eqnarray*}
\int_{\mR^d}|\rho_\varepsilon(t,x)|dx\leq \int_{\mR^d}|\rho_0(x)|dx,
\end{eqnarray*}
which reveals that the first inequality in (\ref{3.27}) is valid.
 \vskip1mm\par
If we set $\rho_\varepsilon^h(t,x)=\rho_\varepsilon(t,x+h)$ and $\rho_0^h(x)=\rho_0(x+h)$, then
\begin{eqnarray*}
\left\{ \begin{array}{ll}
\frac{\partial }{\partial t}\rho_\varepsilon^h(t,x)+\div_xF(\rho_\varepsilon^h)
+(-\Delta_x)^{\frac{\alpha}{2}}\rho_\varepsilon^h-
\varepsilon\Delta\rho_\varepsilon^h=0, \ \ (t,x)\in (0,T) \times \mR^d,
\\ \rho_\varepsilon^h(t=0)=\rho_0^h, \ \ x\in\mR^d.
\end{array}\right.
\end{eqnarray*}
An analogue discussion (as used from (\ref{3.2}) to (\ref{3.24})  leads to
\begin{eqnarray}\label{3.29}
\int_{\mR^d}|\rho_\varepsilon^h(t,x)-\rho_\varepsilon(t,x)|dx\leq
\int_{\mR^d}|\rho_0^h(x)-\rho_0(x)|dx.
\end{eqnarray}
So the second inequality in (\ref{3.27}) satisfies if taking $h$   to zero.
 \vskip1mm\par
The third inequality in (\ref{3.27}) is from the following estimate:
\begin{eqnarray}\label{3.30}
\|\partial_t\rho_\varepsilon(t)\|_{L^1(\mR^d)}
\leq \|\div_xF(\rho_0)
+(-\Delta_x)^{\frac{\alpha}{2}}\rho_0\|_{L^1(\mR^d)} \leq
C([\|\rho_0\|_{L^1(\mR^d)}+\|\rho_0\|_{BV(\mR^d)}].
\end{eqnarray}
 \vskip1mm\par
By (\ref{3.27}), using the Helly theorem (see \cite{Daf} $p_{17}$), the Fr\'{e}chet-Kolmogorov compactness theorem (see \cite{Yos} $p_{275}$) and the Arzela-Ascoli compactness criterion (see \cite{Yos} $p_{85}$), after a standard control of decay at infinity, there is a subsequence (denoted by itself), such that
\begin{eqnarray}\label{3.31}
\rho_\varepsilon\rightarrow \rho \ \mbox{in} \ \cC([0,T];L^1(\mR^d)), \  \mbox{as} \ \varepsilon\rightarrow 0.
\end{eqnarray}
By (\ref{3.31}), from (\ref{3.27}), if we let $\varepsilon\rightarrow 0$, then
\begin{eqnarray*}
\int_{\mR^d}|\rho^h(t,x)-\rho(t,x)|dx\leq
\int_{\mR^d}|\rho_0^h(x)-\rho_0(x)|dx,
\end{eqnarray*}
which implies
\begin{eqnarray}\label{3.32}
\rho\in L^\infty([0,T);BV(\mR^d)).
\end{eqnarray}
 \vskip1mm\par
Besides, $\rho$ satisfy (\ref{2.1})-(\ref{2.3}) for any smooth convex function $\eta$. So $\rho$ is an entropy solution of (\ref{1.1})-(\ref{1.2}). Then Theorem 2.1 applies and thus  $\rho$ is a kinetic solution.
  \vskip1mm\par
For $\rho_0\in L^1\cap BV(\mR^d)$, we approximate it by $\rho_0^\sigma\in L^\infty\cap L^1\cap BV(\mR^d)$, such that
\begin{eqnarray}\label{3.33}
\rho_0^\sigma\rightarrow \rho_0 \ \ \mbox{in} \ L^1\cap BV(\mR^d), \ \ \|\rho_0^{\sigma}\|_{L^1(\mR^d)}\leq C\|\rho_0\|_{L^1(\mR^d)}, \ \|\rho_0^{\sigma}\|_{BV(\mR^d)}\leq C\|\rho_0\|_{BV(\mR^d)}.
\end{eqnarray}
Then there is a kinetic solution $\rho_\sigma$ of (\ref{1.1})-(\ref{1.2}), and for any $h\in\mR$,
\begin{eqnarray}\label{3.34}
\int_{\mR^d}|\rho_{\sigma}^h(t,x)-\rho_{\sigma}(t,x)|dx\leq
\int_{\mR^d}|\rho_0^{\sigma}(x+h)-\rho_0^{\sigma}(x)|dx.
\end{eqnarray}
Correspondingly, the nonnegative measures $m_\sigma$ and $n_\sigma$ meet (\ref{2.14}) and (\ref{2.15}).
 \vskip1mm\par
Moreover for any $\sigma_1, \sigma_2>0$,
\begin{eqnarray}\label{3.35}
\|\rho_{\sigma_1}(t)-\rho_{\sigma_2}(t)\|_{L^1(\mR^{d})}\leq  \|\rho_0^{\sigma_1}-\rho_0^{\sigma_2}\|_{L^1(\mR^d)}.
\end{eqnarray}
 \vskip1mm\par
Denote $\cM_b([0,T)\times\mR^d)$ for  the space of bounded Borel measures over $[0,T)\times\mR^d$, with norm given by the total variation of measures), $n\in L^\infty(\mR_v;L^1([0,T)\times\mR^d)$.
With the aid of (\ref{2.14})-(\ref{2.15}) and (\ref{3.34})-(\ref{3.35}), by choosing a subsequence (not labeled), there are $\rho\in \cC([0,T];L^1(\mR^d))\cap L^\infty([0,T);BV(\mR^d)) $, $m\in L^\infty(\mR_v;\cM_b([0,T)\times\mR^d))$ such that
\begin{eqnarray}\label{3.36}
 \mbox{as} \  \sigma \rightarrow 0, \
\left\{ \begin{array}{ll} \rho_\sigma \rightarrow \rho, \ \mbox{in} \  \cC([0,T];L^1(\mR^d)), \\  m_\sigma\rightarrow m, \ \mbox{in} \  L^\infty_w(\mR_v;\cM_b([0,T]\times\mR^d)), \\ n_\sigma\rightarrow n \ \mbox{in} \  L^\infty_w(\mR_v;L^1([0,T)\times\mR^d)),
\end{array}\right.
\end{eqnarray}
and $n$ fulfills (\ref{2.14}), $m+n\in L^\infty(\mR_v;L^1([0,T)\times\mR^d))$.
 \vskip1mm\par
Moreover, by Remark 2.3, if one takes Kru\u{z}kov entropy $(\rho_\sigma-v)_+$ and $(\rho_\sigma-v)_-$, respectively,  then $m_\sigma+n_\sigma$  satisfies  (\ref{2.17}) and  (\ref{2.18}), respectively. Therefore the nonnegative measures $m$ and $n$ fulfilling (\ref{2.15}), and $\rho$ is a kinetic solution of (\ref{1.1})-(\ref{1.2}).

\begin{remark} (i) $m_{i,\varepsilon}(t,x,v) \ (i=1,2)$ meet $(\ref{3.4})_1$ and the left hand side in $(\ref{3.4})_1$ belongs to $L^1_{loc}$ in variable $v$, so $m_{i,\varepsilon}(t,x,v)+n_{i,\varepsilon}(t,x,v)$ are continuous in $v$. Clearly, $n^1_{i,\varepsilon}(t,x,v)$ and $n^2_{i,\varepsilon}(t,x,v)$ are continuous in $v$. Thus $m_{i,\varepsilon}(t,x,v) (i=1,2)$ are continuous in $v$, which suggests  that $m_{2,\varepsilon}(t,x,0)$ in (\ref{3.12}) and $m_{1,\varepsilon}(t,x,0)$ in (\ref{3.13}) are legitimate.
 \vskip1mm\par
(ii) In our proof, we used the functional (see (\ref{3.7}))
\begin{eqnarray}\label{3.37}
G(t,x,v)&=&|\chi_{\rho_1(t,x)}(v)|+|\chi_{\rho_2(t,x)}(v)|
-2\chi_{\rho_1(t,x)}(v)\chi_{\rho_2(t,x)}(v)
\cr\cr&=&|u_1(t,x,v)|+|u_2(t,x,v)|-2u_1(t,x,v)u_2(t,x,v),
\end{eqnarray}
which was introduced by Perthame (consult to \cite{Per,Per1}) for first order hyperbolic equations. Then this method was extended to the hyperbolic-parabolic equations by Chen and Perthame \cite{CP}, to derive the uniqueness for kinetic solutions. Here our proof follows Chen and Perthame's work, by applying the contraction mapping principle to get the uniqueness of kinetic solutions.
 \end{remark}

\begin{remark} Our existence and uniqueness result can be extended in a routine way to the non-homogeneous problem (\ref{2.19}), if we suppose (\ref{2.20}) and
\begin{eqnarray}\label{3.38}
B^\prime\in L^\infty(\mR), \ \mbox{and} \ \exists \ M_1,M_2\in \mR, \  1_{v>0}A^\prime(v)\leq M_1, \ -M_2\leq1_{v\leq0}A^\prime(v)\leq M_1.
\end{eqnarray}
\vskip1mm\par
Indeed, if one takes functional $G$ as in (\ref{3.37}), by repeating the manipulations from (\ref{3.2}) to (\ref{3.23}),  we end up with
\begin{eqnarray*}
\frac{d}{dt} \int_{\mR^{d+1}}|u_1(t)-u_2(t)|dxdv \leq \int_{\mR^{d+1}}A^\prime(v)|u_1(t)-u_2(t)|dxdv,
\end{eqnarray*}
which demonstrates the uniqueness.
 \vskip1mm\par
For existence part, we   choose $\rho_0\in L^\infty\cap L^1\cap BV(\mR^d)$ first, and consider the approximating problem:
\begin{eqnarray}\label{3.39}
\left\{ \begin{array}{ll}
\frac{\partial }{\partial t}\rho_\varepsilon(t,x)+\div_xF(\rho_\varepsilon)
+(-\Delta_x)^{\frac{\alpha}{2}}B(\rho_\varepsilon)-
\varepsilon\Delta\rho_\varepsilon=A(\rho_\varepsilon), \ \ (t,x)\in (0,T) \times \mR^d,
\\ \rho_\varepsilon(t=0)=\rho_0, \ \ x\in\mR^d.
\end{array}\right.
\end{eqnarray}
We can derive an analogue of (\ref{3.27})
\begin{eqnarray}\label{3.40}
\left\{ \begin{array}{ll} \|\rho_\varepsilon(t)\|_{L^1(\mR^d)}\leq  \exp(M_1t) \|\rho_0\|_{L^1(\mR^d)}, \  \|\rho_\varepsilon(t)\|_{BV(\mR^d)}\leq  \exp(M_1t)\|\rho_0\|_{BV(\mR^d)},
 \\
\|\partial_t\rho_\varepsilon(t)\|_{L^1(\mR^d)}\leq C(t,A,\|\rho_0\|_{L^\infty}) (\|\rho_0\|_{BV(\mR^d)}+\|\rho_0\|_{L^1(\mR^d)}),
\end{array}\right.
\end{eqnarray}
and in view of entropy formulation (\ref{2.3}) (also see Remark 2.3 (ii)), it yields
\begin{eqnarray*}
&&\int_0^T\int_{\mR^d}(m+n)(dt,dx,v)\cr\cr&\leq& \frac{1}{2}\int_0^T\int_{\mR^d}\sgn(\rho-v)_+
A(\rho)dxdt+
\frac{1}{2}\int_{\mR^d}
|(\rho_0(x)-v)_+-(\rho(T,x)-v)_+|
dx
\end{eqnarray*}
and
\begin{eqnarray*}
&&\int_0^T\int_{\mR^d}(m+n)(dt,dx,v)\cr\cr&\leq&
\frac{1}{2}\int_0^T\int_{\mR^d}\sgn(\rho-v)_-
A(\rho)dxdt+\frac{1}{2}\int_{\mR^d}|(\rho_0(x)-v)_--(\rho(T,x)-v)_-|dx.
\end{eqnarray*}
Therefore
\begin{eqnarray*}
&&\int_0^T\int_{\mR^d}(m+n)(dt,dx,v)\cr\cr&\leq& \frac{M_1}{2}1_{v>0}\int_0^T\int_{\mR^d}\sgn(\rho-v)_+
|\rho|dxdt+
\frac{1}{2}1_{v>0} \|(\rho_0(x)-v)_+-(\rho(T,x)-v)_+\|_{L^1(\mR^d)} \cr\cr&&+\frac{M_2}{2}1_{v\leq0}\int_0^T\int_{\mR^d}\sgn(\rho-v)_-
|\rho|dxdt+
\frac{1}{2}1_{v\leq0} \|(\rho_0(x)-v)_--(\rho(T,x)-v)_-\|_{L^1(\mR^d)}.
\end{eqnarray*}
Combining a compactness argument,   we complete the proof for regular initial data.
\vskip1mm\par
Secondly, for $\rho_0\in L^1\cap BV(\mR^d)$, by an approximate discussion, we gain an analogue conclusion of (\ref{3.36}).
\end{remark}
  \vskip1mm\par
With the same verification as in Theorem 3.1, we achieve the following result.
\begin{corollary}
\textbf{(Comparison Principle)} Let (\ref{1.3}) (\ref{2.20}) and (\ref{3.38}) hold and $\rho_{0,1},\rho_{0,2}\in L^1\cap BV(\mR^d)$. Assume that $\rho_1$ and $\rho_2$ are two kinetic solutions of $(\ref{2.19})_1$, to initial values $\rho_{0,1}$ and $\rho_{0,2}$ respectively.  Then
\begin{eqnarray}\label{3.41}
\|\rho_1(t)-\rho_2(t)\|_{L^1(\mR^{d})}\leq \exp(M_1t)\|\rho_{0,1}-\rho_{0,2}\|_{L^1(\mR^d)}.
\end{eqnarray}
Besides, if $\rho_{0,1}\leq \rho_{0,2}$, then
$\rho_1\leq \rho_2$ and in particular, if the initial value is nonnegative, the unique kinetic solution is nonnegative as well.
\end{corollary}

\begin{remark} From above comparison principle (\ref{3.41}), if $M_1<0$ (for example $A(\rho)=M_1\rho$), for any initial data $\rho_0\in L^1\cap BV(\mR^d)$, then the unique kinetic solution $\rho$ of
\begin{eqnarray}\label{3.42}
\frac{\partial}{\partial
t} \rho(t,x)+\div_xF(\rho)+(-\Delta_x)^{\frac{\alpha}{2}}
B(\rho)=M_1\rho, \ \ (t,x)\in (0,T) \times \mR^d,
\end{eqnarray}
converges to zero as $t\rightarrow \infty$,  i.e. $\{0\}$ is the unique global attractor for the solution semigroup.
\end{remark}
 \vskip1mm\par
The restriction conditions on $A$ seem to be strict, but
there are   models, in population dynamics, chemical wave propagation and fluid mechanics, satisfying this assumption. We  now illustrate it by an example.
 \vskip1mm\par
\begin{example}
Consider the following multidimensional   fractional
Burgers-Fisher type equation
\begin{eqnarray}\label{3.43}
\left\{ \begin{array}{ll} \frac{\partial }{\partial
t}\rho(t,x)+\div_x(a\rho^\iota)+
\nu(-\Delta_x)^{\frac{\alpha}{2}}\rho=A(\rho), \ \ (t,x)\in (0,T) \times \mR^d, \\ \rho(t=0,x)=\rho_0(x), \ \ x\in \mR^d,
\end{array}\right.
\end{eqnarray}
where $a\in \mR^d$ is a vector, $\beta\geq0$ $\iota\in \mN$ and
\begin{eqnarray}\label{3.44}
A(\rho)=\left\{ \begin{array}{ll} \beta\rho(1-\rho^k), \ \mbox{when} \ \rho\geq0, \\ \ \ \  0, \ \ \ \mbox{otherwise}.
\end{array}\right.
\end{eqnarray}
 \vskip1mm\par
When $\alpha=2$, $a=0$ and $d=k=1$, it is well known as Fisher equation, proposed by \cite{Fis} in population dynamics, where $\nu>0$ is a diffusion constant, $\beta> 0$ is the linear growth rate. When
$\alpha=2$, $d=1$ and $\iota=k$, it is well known as generalized Burgers-Fisher equation, which is modeled for describing the interaction between reaction mechanisms, convection effects and diffusion transports
\cite{IRR}. And when $\alpha\in (0,2)$, $\beta=0$, it is the generalized fractal/fractional Burgers equation appeared in
continuum mechanics and discussed by \cite{BFW}. The aim of this work is to argue the more general form of the Burger-Fisher and fractal/fractional Burgers equations called generalized fractional Burgers-Fisher type equation in order to show the effectiveness of the current method.
 \end{example}
  \vskip1mm\par
Clearly, $F\in W^{1,\infty}_{loc}(\mR^d)$ and when $k$ is even, (\ref{3.38}) holds with $M_1=\beta$, $M_2=0$. By Remark 3.2 and Corollary 3.1, we have the following result.
\begin{corollary} Let $0\leq\rho_0\in L^1\cap BV(\mR^d)$, $\alpha\in (0,1)$ and $k$ be an even number. Then there is a unique kinetic solution to (\ref{3.43}). Besides, the unique kinetic solution is nonnegative as well.
\end{corollary}

\section{Continuous dependence on   nonlinearities and L\'{e}vy measures}\label{sec4}\setcounter{equation}{0}

This section is devoted to discuss the regularity on $t$ and
the continuous dependence on $f$, $\nu$ and $\alpha$. Since the argument for nonhomogeneous problem is similar, we only concentrate our attention on homogeneous case and our main result is given by:
\begin{theorem} Consider the following Cauchy problems
\begin{eqnarray}\label{4.1}
\left\{ \begin{array}{ll} \frac{\partial}{\partial
t} \rho^1_\alpha(t,x)+\div_xF_1(\rho^1_\alpha)+
\nu_1(-\Delta_x)^{\frac{\alpha}{2}}\rho^1_\alpha=0, \ \ (t,x)\in (0,T) \times \mR^d, \\ \rho^1_\alpha(t=0,x)=\rho_0(x), \ \ x\in \mR^d,
\end{array}\right.
\end{eqnarray}
and
\begin{eqnarray}\label{4.2}
\left\{\begin{array}{ll} \frac{\partial}{\partial
t} \rho^2_\beta(t,x)+\div_xF_2(\rho^2_\beta)+
\nu_2(-\Delta_x)^{\frac{\beta}{2}}\rho^2_\beta=0, \ \ (t,x)\in (0,T) \times \mR^d, \\ \rho^2_\beta(t=0,x)=\rho_0(x), \ \ x\in \mR^d,
\end{array}\right.
\end{eqnarray}
where
\begin{eqnarray}\label{4.3}
\rho_0\in L^1\cap BV(\mR^d), \ F_1,F_2\in W^{1,\infty}_{loc}(\mR;\mR^d) \ \mbox{and} \ F^\prime_1-F^\prime_2\in L^\infty(\mR;\mR^d), \ \alpha,  \beta \in (0,1).
\end{eqnarray}
Let $\rho^1_\alpha$, respectively $\rho^2_\beta$, be the unique kinetic solution to (\ref{4.1}),  respectively to (\ref{4.2}). Then the following claims hold:
\vskip1mm\par
(i) $\rho^1_\alpha$ and $\rho^2_\beta$ are Lipschitz continuous in $t$ in the following sense:  For every  $t,s\in [0,T]$,
\begin{eqnarray}\label{4.4}
\|\rho^1_\alpha(t)-\rho^1_\alpha(s)\|_{L^1(\mR^d)}\leq \|\rho_0\|_{BV(\mR^d)}\|F^\prime_1\|_{L^\infty(\mR)}|t-s|+C
\|\rho_0\|_{L^1(\mR^d)}^{1-\alpha}\|\rho_0\|_{BV(\mR^d)}^{\alpha}|t-s|,
\end{eqnarray}
and
\begin{eqnarray}\label{4.5}
\|\rho^2_\beta(t)-\rho^2_\beta(s)\|_{L^1(\mR^d)}\leq \|\rho_0\|_{BV(\mR^d)}\|F^\prime_2\|_{L^\infty(\mR)}|t-s|+C
\|\rho_0\|_{L^1(\mR^d)}^{1-\beta}\|\rho_0\|_{BV(\mR^d)}^{\beta}|t-s|,
\end{eqnarray}
if $F^\prime_1,F^\prime_2\in L^\infty(\mR;\mR^d)$;
\vskip1mm\par
(ii) Continuous in the nonlinearities and viscosity coefficients: If $\alpha=\beta$, then
\begin{eqnarray}\label{4.6}
\|\rho^1_\alpha-\rho^2_\alpha\|_{\cC([0,T];L^1(\mR^d))}\leq T\|\rho_0\|_{BV(\mR^d)}
\|F^\prime_1-F^\prime_2\|_{L^\infty(\mR)}+T|\nu_1-\nu_2|\|\rho_0\|_{L^1(\mR^d)}^{1-\alpha}
\|\rho_0\|_{BV(\mR^d)}^{\alpha};
\end{eqnarray}
\vskip1mm\par
(iii) Lipschitz continuous in L\'{e}vy measure: If $F_1=F_2$, then for every $\lambda\in(0,1)$
\begin{eqnarray}\label{4.7}
\limsup_{\alpha,\beta\rightarrow\lambda}\frac{\|\rho_\alpha^1-
\rho_\beta^1\|_{\cC([0,T];L^1(\mR^d))}}{|\alpha-\beta|}
\leq CT\|\rho_0\|_{L^1(\mR^d)}^{1-\lambda}\|\rho_0\|_{BV(\mR^d)}^{\lambda}
(1+|\log\frac{\|\rho_0\|_{L^1}}{\|\rho_0\|_{BV}}|).
\end{eqnarray}
\end{theorem}
 \vskip1mm\par
Before proving the main result,   we introduce another notion of   solutions and present a useful lemma.
\begin{definition} Let $\rho_0\in L^\infty(\mR^d), \ F_1\in W^{1,\infty}_{loc}(\mR;\mR^d)$ and $
\alpha\in (0,1)$. We call $\rho\in L^\infty([0,T)\times \mR^d)$ an entropy solution of (\ref{4.1}), if for every $v\in\mR, r>0$ and every  nonnegative function $\psi\in\cD([0,T)\times\mR^d)$
\begin{eqnarray}\label{4.8}
&&\int_0^T\int_{\mR^d}{\Big [ }|\rho^1_\alpha-v| \frac{\partial}{\partial t}\psi(t,x)+\sgn(\rho^1_\alpha-v)(F_1(\rho^1_\alpha)-F_1(v))\cdot\nabla\psi{\Big ]}dtdx +\int_{\mR^d}|\rho_0-v|\psi(0,x)dx\cr\cr&&+\int_0^T\int_{\mR^d}{\Big [ }|\rho^1_\alpha-v| \cL^\alpha_r[\psi(t,x)]+\sgn(\rho^1_\alpha-v)\cL^{\alpha,r}[\rho^1_\alpha]\psi(t,x){\Big ]}dtdx\geq0,
\end{eqnarray}
where $\cL^\alpha_r$ and $\cL^{\alpha,r}$ are defined, for   $\varphi \in \cD(\mR^d)$ and $x\in \mR^d$, by
\begin{eqnarray*}
&&\cL^\alpha_r\varphi(x)=
c(d,\alpha)\int_{|z|<r}\frac{\varphi(x+z)-\varphi(x)-\nabla\varphi(x)\cdot z}
{|z|^{d+\alpha}}dz,
\cr\cr&& \quad\quad
\cL^{\alpha,r}\varphi(x)=
c(d,\alpha)\int_{|z|<r}\frac{\varphi(x+z)-\varphi(x)}
{|z|^{d+\alpha}}dz.
\end{eqnarray*}
\end{definition}

\begin{lemma} (i) (\cite{CJ} Theorem 2.5) If $\rho_0\in L^1\cap L^\infty\cap BV(\mR^d)$, then   Definition 2.1 and Definition 4.1 are equivalent. Thus by the kinetic formulation in Theorem 2.1, both  Definition 2.2 and Definition 4.1 are equivalent.
 \vskip1mm\par
(ii) (\cite{ACJ1} Theorem 3.3, Theorem 3.4 or \cite{ACJ2} Theorem 2, Theorem 4) If $\rho_0\in L^1\cap L^\infty\cap BV(\mR^d)$, then Theorem 4.1 holds.
\end{lemma}

\begin{remark} Notice that the right hand sides in (\ref{4.4})-(\ref{4.7}) are dependent only on $\|\rho_0\|_{L^1}$ and $ \|\rho_0\|_{BV}$. So (\ref{4.4})-(\ref{4.7})  may be true if the entropy solutions take values in $\cC([0,T];L^1(\mR^d))\cap L^\infty([0,T);BV(\mR^d))$. However, when $\rho\in \cC([0,T];L^1(\mR^d))\cap L^\infty([0,T);BV(\mR^d))$, the term $ \int_0^T\int_{\mR^d}\sgn(\rho^1_\alpha-v)(F_1(\rho^1_\alpha)-F_1(v))\cdot\nabla\psi(t,x) dtdx$ in the first line in (\ref{4.7}) is not legitimate. To overcome this obstacle,  we introduce the notion of kinetic solutions, and
extend Lemma 4.1 to the class of $L^1\cap BV(\mR^d)$ solutions.
\end{remark}
 \vskip0mm\noindent
\textbf{Proof of Theorem 4.1.}  We approximate $\rho_0$ by $\rho_0^\sigma$ such that (\ref{3.30}) holds. By Remark 2.3, (\ref{2.9}) and (\ref{3.30}), we end up with
\begin{eqnarray}\label{4.9}
\|\rho^{1,\sigma}_\alpha(t)-\rho^{1,\sigma}_\alpha(s)\|_{L^1(\mR^d)}\leq \|\rho_0^{\sigma}\|_{BV(\mR^d)}\|F^\prime_1\|_{L^\infty(\mR)}|t-s|+C
\|\rho_0^{\sigma}\|_{L^1(\mR^d)}^{1-\alpha}\|\rho_0^{\sigma}\|_{BV(\mR^d)}^{\alpha}|t-s|
\end{eqnarray}
and
\begin{eqnarray}\label{4.10}
\|\rho^{2,\sigma}_\beta(t)-\rho^{2,\sigma}_\beta(s)\|_{L^1(\mR^d)}\leq \|\rho_0^{\sigma}\|_{BV(\mR^d)}\|F^\prime_2\|_{L^\infty(\mR)}|t-s|+C
\|\rho_0^{\sigma}\|_{L^1(\mR^d)}^{1-\beta}\|\rho_0^{\sigma}\|_{BV(\mR^d)}^{\beta}|t-s|,
\end{eqnarray}
where the constant $C$ is dependent only on $\|F^\prime_1\|_{L^\infty(\mR)}, \|F^\prime_2\|_{L^\infty(\mR)}, d$ and $\alpha,\beta$.
 \vskip1mm\par
By Lemma 4.1 (i) and and (ii), $\rho^{1,\sigma}_\alpha$ and $\rho^{2,\sigma}_\beta$ fulfill
\begin{eqnarray}\label{4.11}
\|\rho^{1,\sigma}_\alpha\!-\!\rho^{2,\sigma}_\alpha\|_{\cC([0,T];L^1(\mR^d))}\!\leq\! T\|\rho_0^{\sigma}\|_{BV(\mR^d)}
\|F^\prime_1-F^\prime_2\|_{L^\infty(\mR)}\!+\!T|\nu_1-\nu_2|
\|\rho_0^\sigma\|_{L^1(\mR^d)}^{1-\alpha}
\|\rho_0^\sigma\|_{BV(\mR^d)}^{\alpha}
\end{eqnarray}
and
\begin{eqnarray}\label{4.12}
\|\rho^{1,\sigma}_\alpha-
\rho_\beta^{1,\sigma}\|_{\cC([0,T];L^1(\mR^d))}
\leq T\int_{\mR^d}\|\rho_0^{\sigma}(\cdot+z)-\rho_0^{\sigma}(\cdot)\|_{L^1(\mR^d)}
d|\mu_\alpha-\mu_\beta|,
\end{eqnarray}
where
$$
d\mu_\alpha=\frac{c(d,\alpha)}{|z|^{d+\alpha}}dz, \ d\mu_\beta=\frac{c(d,\beta)}{|z|^{d+\beta}}dz.
$$
 \vskip1mm\par
Observing that
$$
\rho^{1,\sigma}_\alpha  \rightarrow \rho^1_\alpha,
\rho^{2,\sigma}_\beta \rightarrow \rho^2_\beta, \ \rho^{2,\sigma}_\alpha  \rightarrow \rho^2_\alpha, \ \rho^{1,\sigma}_\beta  \rightarrow \rho^1_\beta \  \mbox{in} \  \cC([0,T];L^1(\mR^d)), \  \mbox{as} \  \sigma \rightarrow 0,
$$
as $\sigma \rightarrow 0$, and noting  (\ref{4.9})-(\ref{4.12}), we arrive at inequalities (\ref{4.4})-(\ref{4.6}). Therefore, the claims $(i)$ and $(ii)$ in Theorem 4.1 hold.
 \vskip1mm\par
From (\ref{4.9})-(\ref{4.12}), we also have
\begin{eqnarray}\label{4.13}
\|\rho^1_\alpha-
\rho_\beta^1\|_{\cC([0,T];L^1(\mR^d))}
\leq T\int_{\mR^d}\|\rho_0(\cdot+z)-\rho_0(\cdot)\|_{L^1(\mR^d)}
d|\mu_\alpha-\mu_\beta|,
\end{eqnarray}
and to prove claim $(iii)$, let $0<r_1\in \mR$, we split the integral in the right hand side in (\ref{4.13}) into two parts
\begin{eqnarray}\label{4.14}
\int_{|z|\geq r_1}\|\rho_0(\cdot+z)-\rho_0(\cdot)\|_{L^1(\mR^d)}
d|\mu_\alpha-\mu_\beta|+\int_{|z|< r_1}\|\rho_0(\cdot+z)-\rho_0(\cdot)\|_{L^1(\mR^d)}
d|\mu_\alpha-\mu_\beta|.
\end{eqnarray}
Then the proof for Theorem 4 (\cite{ACJ2}) applies, and we obtain (\ref{4.6}).  This completes the proof.

\begin{remark} We can prove the continuous dependence of solutions on nonlinearities by introducing the functional $G$ (given in (\ref{3.37})). Indeed, we write (\ref{4.1}) and (\ref{4.2}) in microscopic types by using kinetic formulation first, then we regularize solutions in $t,x,v$ and repeat the calculations from (\ref{3.11}) to (\ref{3.23}), to get
\begin{eqnarray}\label{4.15}
&&\frac{d}{dt}\int_{\mR^{d+1}}|u_\alpha^1(t,x,v)-u_\alpha^2(t,x,v)|dxdv\cr\cr&\leq& 2\|F^\prime_1-F^\prime_2\|_{L^\infty(\mR)}\|\rho_0\|_{BV(\mR^d)}+2|\nu_1-\nu_2|\|\rho_0\|_{L^1(\mR^d)}^{1-\alpha}
\|\rho_0\|_{BV(\mR^d)}^{\alpha},
\end{eqnarray}
where $u_\alpha^1=\chi_{\rho_\alpha^1}(v)$ and $u_\alpha^2=\chi_{\rho_\alpha^2}(v)$.
 \vskip1mm\par
From (\ref{4.15}), it follows that
\begin{eqnarray}\label{4.16}
\|\rho^1_\alpha-\rho^2_\alpha\|_{\cC([0,T];L^1(\mR^d))}\!\leq\! 2T[\|\rho_0\|_{BV(\mR^d)}
\|F^\prime_1-F^\prime_2\|_{L^\infty(\mR)}\!+\!|\nu_1-\nu_2|\|\rho_0\|_{L^1(\mR^d)}^{1-\alpha}
\|\rho_0\|_{BV(\mR^d)}^{\alpha}].
\end{eqnarray}
If we define the right hand side of (\ref{4.16}) by $I(T,\rho_0,F_1,F_2,\nu_1,\nu_2)$, then from (\ref{4.6}),
$$
\|\rho^1_\alpha-\rho^2_\alpha\|_{\cC([0,T];L^1(\mR^d))}
\leq\frac{1}{2}I(T,\rho_0,F_1,F_2,\nu_1,\nu_2).
$$
So (\ref{4.6}) implies (\ref{4.16}), and in this sense, we say the estimate (\ref{4.6}) is better than (\ref{4.16}). Hence in the proof of Theorem 4.1, we adapt the method developed in \cite{ACJ1,ACJ2}.
\end{remark}
\vskip1mm\par
Besides the continuous dependence, we also have obtained the limiting equations as $\alpha\downarrow0$ and $\nu\downarrow0$. Firstly, we give a useful lemma for fixed $\nu$, which will serve us well for the limiting problem as $\alpha\downarrow0$, and for simplicity we take $\nu=1$.
\begin{lemma} (\cite{ACJ2} Theorem 3) Let $\rho_0\in L^\infty(\mR^d)$ and for $\alpha\in(0,1)$, let $\rho_\alpha$ be the unique entropy solution (defined by Definition 4.1) of (\ref{1.1})-(\ref{1.2}).
If $\rho_0\in L^1(\mR^d)$, then as $\alpha\downarrow0$, $\rho_\alpha$ converges in $\cC([0,T];L^1_{loc}(\mR^d))$ to the unique entropy solution (defined by Definition 4.1) $\rho\in L^\infty([0,T)\times \mR^d)\cap \cC([0,T];L^1(\mR^d))$ of the Cauchy problem
\begin{eqnarray}\label{4.17}
\left\{ \begin{array}{ll} \frac{\partial}{\partial
t} \rho(t,x)+\div_xF(\rho)+
\rho=0, \ \ (t,x)\in (0,T) \times \mR^d, \\ \rho(t=0,x)=\rho_0(x), \ \ x\in \mR^d.
\end{array}\right.
\end{eqnarray}
\end{lemma}
 \vskip1mm\par
Our main result is given by:
\begin{theorem}
Let $\rho_0\in L^1\cap BV(\mR^d)$,  and for $\alpha\in(0,1)$, let $\rho_\alpha^\nu$ be the unique kinetic solution  of (\ref{1.1})- (\ref{1.2}).
 \vskip1mm\par
(i) As $\nu\downarrow0$, $\rho_\alpha^\nu$ converges in $\cC([0,T];L^1(\mR^d))$ to the unique kinetic solution $\rho$ of
the following Cauchy problem
\begin{eqnarray}\label{4.18}
\left\{ \begin{array}{ll} \frac{\partial}{\partial
t} \rho(t,x)+\div_xF(\rho)=0, \ \ (t,x)\in (0,T) \times \mR^d, \\ \rho(t=0,x)=\rho_0(x), \ \ x\in \mR^d.
\end{array}\right.
\end{eqnarray}
Moreover, we have the following error estimate: for all $T > 0$,
\begin{eqnarray}\label{4.19}
\|\rho_\alpha^\nu-\rho\|_{\cC([0,T];L^1(\mR^d))}=O(\nu), \ as \ \nu\rightarrow 0.
\end{eqnarray}
  \vskip1mm\par
(ii)  If $\alpha\downarrow0$, then $\rho_{\alpha}^\nu$ converges in $\cC([0,T];L^1_{loc}(\mR^d))$ to the unique kinetic solution $\rho^\nu$ of the following Cauchy problem
\begin{eqnarray}\label{4.20}
\left\{ \begin{array}{ll} \frac{\partial}{\partial
t} \rho^\nu(t,x)+\div_xF(\rho^\nu)+
\nu\rho^\nu=0, \ \ (t,x)\in (0,T) \times \mR^d, \\ \rho^\nu(t=0,x)=\rho_0(x), \ \ x\in \mR^d.
\end{array}\right.
\end{eqnarray}
\end{theorem}
\textbf{Proof.} For every pair of $\nu_1,\nu_2>0$, by virtue of Theorem 4.1 (ii), we have
\begin{eqnarray}\label{4.21}
\|\rho_\alpha^{\nu_1}-\rho_\alpha^{\nu_2}\|_{\cC([0,T];L^1(\mR^d))}\leq T|\nu_1-\nu_2|\|\rho_0\|_{L^1(\mR^d)}^{1-\alpha}
\|\rho_0\|_{BV(\mR^d)}^{\alpha},
\end{eqnarray}
which implies that $\{\rho_\alpha^\nu\}_\nu$ is a Cauchy sequence in $\cC([0,T];L^1(\mR^d))$. So $\{u_\alpha^\nu=\chi_{\rho_\alpha^\nu}\}_\nu$ is a Cauchy sequence in $\cC([0,T];L^1(\mR^{d+1}))$.
 \vskip1mm\par
Observe that $u_\alpha^\nu$ yields
\begin{eqnarray}\label{4.22}
\left\{ \begin{array}{ll} \frac{\partial}{\partial t}u_\alpha^\nu(t,x,v)+f(v)\cdot \nabla_xu_\alpha^\nu+ \nu(-\Delta_x)^{\frac{\alpha}{2}}u_\alpha^\nu=\frac{\partial}{\partial v}(m_\alpha^\nu+n_\alpha^\nu), \ (t,x,v)\in (0,T) \times \mR^d\times \mR,  \\ u_\alpha^\nu(t=0,x,v)=\chi_{\rho_0(x)}(v), \ \ x\in \mR^d.
\end{array}\right.
\end{eqnarray}
Combining (\ref{2.14}) and (\ref{2.17}), we conclude that
\begin{eqnarray}\label{4.23}
n_\alpha^\nu\rightarrow 0 \ \mbox{in} \ L^1_{loc}(\mR_v;L^1([0,T]\times\mR^d)), \ \mbox{as} \ \nu\downarrow0.
\end{eqnarray}
\vskip1mm\par
In view of (\ref{2.15}), (\ref{2.17}) and (\ref{2.18}), there is a nonnegative measure $m\in L^\infty_0(\mR_v;\cM_b([0,T]\times\mR^d))$, so that
\begin{eqnarray}\label{4.24}
m_\alpha^\nu\rightarrow m, \ \mbox{in} \  L^\infty_w(\mR_v;\cM_b([0,T]\times\mR^d)), \ \mbox{as} \ \nu\downarrow0.
\end{eqnarray}
By (\ref{4.21}), (\ref{4.23}), (\ref{4.24}) and the following estimate
\begin{eqnarray}\label{4.25}
\|\rho_\alpha^\nu(t)\|_{BV(\mR^d)}\leq \|\rho_0\|_{BV(\mR^d)},
\end{eqnarray}
and take $\nu\downarrow0$ in (\ref{4.22}) in the distributions sense, we know that there is $\rho\in \cC([0,T];L^1(\mR^d))\cap L^\infty([0,T);BV(\mR^d))$, satisfying
\begin{eqnarray}\label{4.26}
\left\{ \begin{array}{ll} \frac{\partial}{\partial t}\chi_{\rho(t,x)}(v)+f(v)\cdot \nabla_x\chi_{\rho(t,x)}(v)=\frac{\partial}{\partial v}m(t,x,v), \ (t,x,v)\in (0,T) \times \mR^d\times \mR,  \\ u(t=0,x,v)=\chi_{\rho_0(x)}(v), \ \ x\in \mR^d.
\end{array}\right.
\end{eqnarray}
Clearly the kinetic solution for (\ref{4.18}) is unique, and thus $\rho$ is the unique kinetic solution of (\ref{4.18}).
 \vskip1mm\par
The error estimate (\ref{4.19}) follows from (\ref{4.21}) by letting $\nu_2\downarrow0$ and replacing $\nu_1$ by $\nu$,   and this finishes the proof for (i).
 \vskip1mm\par
It remains to show (ii) and without loss of generality, we suppose $\nu=1$.
 \vskip1mm\par
Let $\rho_0^\sigma$ and $\rho^{\sigma}_\alpha$ be described in (\ref{4.9}). Then, by Lemma 4.2, as $\alpha\downarrow0$, $\rho^{\sigma}_\alpha$ converges in $\cC([0,T];L^1_{loc}(\mR^d))$ to the unique entropy solution (defined by Definition 4.1) $\rho^{\sigma}\in L^\infty([0,T)\times \mR^d)\cap \cC([0,T];L^1(\mR^d))$ of
\begin{eqnarray}\label{4.27}
\left\{ \begin{array}{ll} \frac{\partial}{\partial
t} \rho^{\sigma}(t,x)+\div_xF(\rho^{\sigma})+
\rho^{\sigma}=0, \ \ (t,x)\in (0,T) \times \mR^d, \\ \rho^{\sigma}(t=0,x)=\rho_0^\sigma(x), \ \ x\in \mR^d,
\end{array}\right.
\end{eqnarray}
 \vskip1mm\par
With the aid of classical kinetic formulation (see \cite{Per1}), $\rho^{\sigma}\in L^\infty([0,T)\times \mR^d)\cap \cC([0,T];L^1(\mR^d))\cap L^\infty([0,T);BV(\mR^d))$ and it is the unique kinetic solution of (\ref{4.27}), i.e. $u^\sigma(t,x,v)=\chi_{\rho^{\sigma}}(v)$ meets
\begin{eqnarray}\label{4.28}
\left\{ \begin{array}{ll} \frac{\partial}{\partial t}u^\sigma+f(v)\cdot \nabla_xu^\sigma-v\frac{\partial}{\partial v}u^\sigma=\frac{\partial}{\partial v}m^\sigma, \ (t,x,v)\in (0,T) \times \mR^d\times \mR, \\ u^\sigma(t=0)=\chi_{\rho_0^\sigma}(v), \  (x,v)\in \mR^d\times \mR, \end{array}\right.
\end{eqnarray}
for some nonnegative measure $m^\sigma$, which satisfies
\begin{eqnarray}\label{4.29}
\int_0^T\int_{\mR^d}m^\sigma(dt,dx,v)\in L^\infty_0(\mR).
\end{eqnarray}
\vskip1mm\par
In view of (\ref{2.15}), (\ref{2.17}) and (\ref{2.18}), there is a nonnegative measure $m\in L^\infty_0(\mR_v;\cM_b([0,T]\times\mR^d))$, so that
\begin{eqnarray}\label{4.30}
m^\sigma\rightarrow m, \ \mbox{in} \  L^\infty_w(\mR_v;\cM_b([0,T]\times\mR^d)), \ \mbox{as} \ \sigma\downarrow0.
\end{eqnarray}
By (\ref{4.24}), (\ref{4.29})-(\ref{4.30}), if we take $\sigma\downarrow0$ in (\ref{4.28}) in the distributions sense, then there is $\rho\in \cC([0,T];L^1(\mR^d))\cap L^\infty([0,T);BV(\mR^d))$, satisfying
\begin{eqnarray}\label{4.31}
\left\{ \begin{array}{ll} \frac{\partial}{\partial t}\chi_{\rho(t,x)}(v)+f(v)\cdot \nabla_x\chi_{\rho(t,x)}(v)-v\frac{\partial}{\partial v}\chi_{\rho}(v)=\frac{\partial}{\partial v}m(t,x,v), \ (t,x,v)\in (0,T) \times \mR^d\times \mR,  \\ u(t=0,x,v)=\chi_{\rho_0(x)}(v), \ \ x\in \mR^d.
\end{array}\right.
\end{eqnarray}
Thus $\rho$ is the unique kinetic solution of (\ref{4.26}) and we complete the proof.

\begin{remark} The calculations for Corollary 3.1 used here, we gain: if $\rho_0\geq 0$, then the unique kinetic solution $\rho$ for (\ref{4.18}), and the unique kinetic solution $\rho^\nu$ for (\ref{4.20})
are nonnegative. Besides, we have the following identities
\begin{eqnarray}\label{4.31}
\int_{\mR^d}\rho(t,x)dx=\int_{\mR^d}\rho_0(x)dx, \ \ \int_{\mR^d}\rho^\nu(t,x)dx+\nu\int_0^t\int_{\mR^d}\rho^\nu(s,x)dxds=
\int_{\mR^d}\rho_0(x)dx.
\end{eqnarray}
On the other hand, if we let $\rho_\alpha^\nu$ be the unique kinetic solution of
\begin{eqnarray}\label{4.33}
\left\{ \begin{array}{ll} \frac{\partial}{\partial
t} \rho^\nu_\alpha(t,x)+\div_xF(\rho^\nu_\alpha)+
\nu(-\Delta_x)^{\frac{\alpha}{2}}\rho^\nu_\alpha=0, \ \ (t,x)\in (0,T) \times \mR^d, \\ \rho^\nu_\alpha(t=0,x)=\rho_0(x)\geq 0, \ \ x\in \mR^d,
\end{array}\right.
\end{eqnarray}
then
\begin{eqnarray}\label{4.31}
\int_{\mR^d}\rho^\nu_\alpha(t,x)dx=\int_{\mR^d}\rho_0(x)dx.
\end{eqnarray}
Therefore, the mass preserving property still holds at the $\nu\downarrow0$ limit, but will be lost at the $\alpha\downarrow0$ limit.
So, in general speaking, as $\alpha\downarrow0$, $\rho_{\alpha}^\nu$ does not converges in $\cC([0,T];L^1(\mR^d))$ to the unique kinetic solution $\rho^\nu$ of (\ref{4.20}). From this point, the convergence here is sharp.  But when discussing (i), the mass preserving property still holds at the limit, so one can expect $L^1$ convergence for (i) as $\nu\downarrow0$. Moreover, the preceding convergence is    in $L^1$ spaces,  but   $L^1\cap BV$ is a proper space to ensure this discussion. Based upon this point, we derive analogue results of Theorem 3.3 \cite{Ali} and Theorem 3 \cite{ACJ2} for kinetic solutions, without assuming $\rho_0\in L^\infty$.
\end{remark}

\vskip4mm\noindent
\textbf{\large{Acknowledgements}}
 \vskip3mm\par
This research was partly supported by the NSF of China grants 11501577, 11301146, 11531006, 11371367  and   11271290.

\end{document}